\documentclass[11pt]{article}
 
 \usepackage[dvips]{graphicx} 
 \usepackage{textcomp}
\usepackage{amsbsy}
\usepackage{latexsym}
\usepackage[mathscr]{eucal}
\usepackage{amsfonts}
\usepackage{amssymb}
\usepackage[usenames]{color}
 \usepackage[dvipsnames]{xcolor}

\usepackage{fullpage}

\begin{document}
 
\bibliographystyle{plain}
 
\title
{
  Data Assimilation in Reduced Modeling
}
\author{ 
Peter Binev, Albert Cohen, Wolfgang Dahmen,\\
 Ronald DeVore, Guergana Petrova, and Przemyslaw Wojtaszczyk
\thanks{
   This research was supported by the ONR Contracts
  N00014-11-1-0712,  N00014-12-1-0561, N00014-15-1-2181; the  NSF Grants  DMS 1222715,  
       DMS 1222390;
  the Institut Universitaire de France; the ERC Adv grant BREAD; 
    the DFG SFB-Transregio 40; 
    the DFG Research Group 1779;
   the Excellence Initiative of the German Federal and State Governments,
   (RWTH Aachen  Distinguished Professorship, Graduate School AICES);
  and the  Polish NCN grant DEC2011/03/B/ST1/04902. }  }

\hbadness=10000
\vbadness=10000

\newtheorem{lemma}{Lemma}[section]
\newtheorem{prop}[lemma]{Proposition}
\newtheorem{cor}[lemma]{Corollary}
\newtheorem{theorem}[lemma]{Theorem}
\newtheorem{remark}[lemma]{Remark}
\newtheorem{example}[lemma]{Example}
\newtheorem{definition}[lemma]{Definition}
\newtheorem{proper}[lemma]{Properties}
\newtheorem{assumption}[lemma]{Assumption}

 \newenvironment{disarray}{\everymath{\displaystyle\everymath{}}\array}{\endarray}

\def\RR{\rm \hbox{I\kern-.2em\hbox{R}}}
\def\NN{\rm \hbox{I\kern-.2em\hbox{N}}}
\def\ZZ{\rm {{\rm Z}\kern-.28em{\rm Z}}}
\def\CC{\rm \hbox{C\kern -.5em {\raise .32ex \hbox{$\scriptscriptstyle
|$}}\kern
-.22em{\raise .6ex \hbox{$\scriptscriptstyle |$}}\kern .4em}}
\def\vp{\varphi}
\def\<{\langle}
\def\>{\rangle}
\def\t{\tilde}
\def\i{\infty}
\def\e{\varepsilon}
\def\sm{\setminus}
\def\nl{\newline}
\def\o{\overline}
\def\wt{\widetilde}
\def\wh{\widehat}
\def\cT{{\cal T}}
\def\cA{{\cal A}}
\def\cI{{\cal I}}
\def\cV{{\cal V}}
\def\cB{{\cal B}}
\def\cF{{\cal F}}

\def\cR{{\cal R}}
\def\cD{{\cal D}}
\def\cP{{\cal P}}
\def\cJ{{\cal J}}
\def\cM{{\cal M}}
\def\cO{{\cal O}}
\def\Chi{\raise .3ex
\hbox{\large $\chi$}} \def\vp{\varphi}
\def\lsima{\hbox{\kern -.6em\raisebox{-1ex}{$~\stackrel{\textstyle<}{\sim}~$}}\kern -.4em}
\def\lsim{\hbox{\kern -.2em\raisebox{-1ex}{$~\stackrel{\textstyle<}{\sim}~$}}\kern -.2em}
\def\[{\Bigl [}
\def\]{\Bigr ]}
\def\({\Bigl (}
\def\){\Bigr )}
\def\[{\Bigl [}
\def\]{\Bigr ]}
\def\({\Bigl (}
\def\){\Bigr )}
\def\L{\pounds}
\def\pr{{\rm Prob}}
 
\newcommand{\cs}[1]{{\color{magenta}{#1}}}
\def\ds{\displaystyle}
\def\ev#1{\vec{#1}}      
\newcommand{\lt}{\ell^{2}(\nabla)}
\def\Supp#1{{\rm supp\,}{#1}}
 
\def\R{\mathbb{R}}
 
\def\E{\mathbb{E}}
\def\nl{\newline}
\def\T{{\relax\ifmmode I\!\!\hspace{-1pt}T\else$I\!\!\hspace{-1pt}T$\fi}}
 
\def\N{\mathbb{N}}
 
\def\Z{\mathbb{Z}}
\def\N{\mathbb{N}}
\def\Zd{\Z^d}
 
\def\Q{\mathbb{Q}}
 \def\C{\mathbb{C}}
\def\Rd{\R^d}
\def\gsim{\mathrel{\raisebox{-4pt}{$\stackrel{\textstyle>}{\sim}$}}}
\def\sime{\raisebox{0ex}{$~\stackrel{\textstyle\sim}{=}~$}}
\def\lsim{\raisebox{-1ex}{$~\stackrel{\textstyle<}{\sim}~$}}
\def\div{\mbox{ div }}
\def\M{M}  \def\NN{N}                  
\def\L{{\ell}}                
\def\Le{{\ell^1}}             
\def\Lz{{\ell^2}}
\def\Let{{\tilde\ell^1}}      
\def\Lzt{{\tilde\ell^2}}
\def\Ltw{\ell^\tau^w(\nabla)}
\def\t#1{\tilde{#1}}
\def\la{\lambda}
\def\La{\Lambda}
\def\ga{\gamma}
\def\BV{{\rm BV}}
\def\Ga{\eta}
\def\al{\alpha}
\def\cZ{{\cal Z}}
\def\cA{{\cal A}}
\def\cU{{\cal U}}
\def\argmin{\mathop{\rm argmin}}
\def\argmax{\mathop{\rm argmax}}
\def\prob{\mathop{\rm prob}}
\def\A{\mathop{\rm Alg}}

\def \bphi{{\bf\phi}}

\def\cO{{\cal O}}
\def\cA{{\cal A}}
\def\cC{{\cal C}}
\def\cS{{\cal F}}
\def\bu{{\bf u}}
\def\bz{{\bf z}}
\def\bZ{{\bf Z}}
\def\bI{{\bf I}}
\def\cE{{\cal E}}
\def\cD{{\cal D}}
\def\cG{{\cal G}}
\def\cI{{\cal I}}
\def\cJ{{\cal J}}
\def\cM{{\cal M}}
\def\cN{{\cal N}}
\def\cT{{\cal T}}
\def\cU{{\cal U}}
\def\cV{{\cal V}}
\def\cW{{\cal W}}
\def\cL{{\cal L}}
\def\cB{{\cal B}}
\def\cG{{\cal G}}
\def\cK{{\cal K}}
\def\cS{{\cal S}}
\def\cP{{\cal P}}
\def\cQ{{\cal Q}}
\def\cR{{\cal R}}
\def\cU{{\cal U}}
\def\bL{{\bf L}}
\def\bl{{\bf l}}
\def\bK{{\bf K}}
\def\bC{{\bf C}}
\def\X{X\in\{L,R\}}
\def\ph{{\varphi}}
\def\D{{\Delta}}
\def\H{{\cal H}}
\def\bM{{\bf M}}
\def\bx{{\bf x}}
\def\bj{{\bf j}}
\def\bG{{\bf G}}
\def\bP{{\bf P}}
\def\bW{{\bf W}}
\def\bT{{\bf T}}
\def\bV{{\bf V}}
\def\bv{{\bf v}}
\def\bt{{\bf t}}
\def\bz{{\bf z}}
\def\bw{{\bf w}}
\def \span{{\rm span}}
\def \meas {{\rm meas}}
\def\rhom{{\rho^m}}
\def\diff{\hbox{\tiny $\Delta$}}
\def\EE{{\rm Exp}}
 
\def\lll{\langle}
\def\argmin{\mathop{\rm argmin}}
\def\argmax{\mathop{\rm argmax}}
 \def\dJ{\nabla}
\newcommand{\ba}{{\bf a}}
\newcommand{\bb}{{\bf b}}
\newcommand{\bc}{{\bf c}}
\newcommand{\bd}{{\bf d}}
\newcommand{\bs}{{\bf s}}
\newcommand{\bff}{{\bf f}}
\newcommand{\bp}{{\bf p}}
\newcommand{\bg}{{\bf g}}
\newcommand{\by}{{\bf y}}
\newcommand{\br}{{\bf r}}
\newcommand{\be}{\begin{equation}}
\newcommand{\ee}{\end{equation}}
\newcommand{\bea}{$$ \begin{array}{lll}}
\newcommand{\eea}{\end{array} $$}
 
\def \Vol{\mathop{\rm  Vol}}
\def \mes{\mathop{\rm mes}}
\def\rad{\mathop{\rm rad}}
\def \Prob{\mathop{\rm  Prob}}
\def \exp{\mathop{\rm    exp}}
\def \sign{\mathop{\rm   sign}}
\newcommand{\mult}{\mathop{\rm   mult}}
\newcommand{\one}{\mathop{\rm   one}}

\def \sp{\mathop{\rm   span}}
\def \vphi{{\varphi}}
\def \csp{\overline \mathop{\rm   span}}
 
\newcommand{\beqn}{\begin{equation}}
\newcommand{\eeqn}{\end{equation}}
\def\beginproof{\noindent{\bf Proof:}~ }
\def\endproof{\hfill\rule{1.5mm}{1.5mm}\\[2mm]}

\newcommand{\utr}{u^{\rm true}}
\newcommand{\Cor}{\kappa}

\newenvironment{Proof}{\noindent{\bf Proof:}\quad}{\endproof}

\renewcommand{\theequation}{\thesection.\arabic{equation}}
\renewcommand{\thefigure}{\thesection.\arabic{figure}}

\makeatletter
\@addtoreset{equation}{section}
\makeatother

\newcommand\abs[1]{\left|#1\right|}
\newcommand\clos{\mathop{\rm clos}\nolimits}
\newcommand\trunc{\mathop{\rm trunc}\nolimits}
\renewcommand\d{d}
\newcommand\dd{d}
\newcommand\diag{\mathop{\rm diag}}
\newcommand\dist{\mathop{\rm dist}}
\newcommand\diam{\mathop{\rm diam}}
\newcommand\cond{\mathop{\rm cond}\nolimits}
\newcommand\eref[1]{{\rm (\ref{#1})}}
\newcommand{\iref}[1]{{\rm (\ref{#1})}}
\newcommand\Hnorm[1]{\norm{#1}_{H^s([0,1])}}
\def\int{\intop\limits}
\renewcommand\labelenumi{(\roman{enumi})}
\newcommand\lnorm[1]{\norm{#1}_{\ell^2(\Z)}}
\newcommand\Lnorm[1]{\norm{#1}_{L_2([0,1])}}
\newcommand\LR{{L_2(\R)}}
\newcommand\LRnorm[1]{\norm{#1}_\LR}
\newcommand\Matrix[2]{\hphantom{#1}_#2#1}
\newcommand\norm[1]{\left\|#1\right\|}
\newcommand\ogauss[1]{\left\lceil#1\right\rceil}
\newcommand{\QED}{\hfill
\raisebox{-2pt}{\rule{5.6pt}{8pt}\rule{4pt}{0pt}}%
  \smallskip\par}
\newcommand\Rscalar[1]{\scalar{#1}_\R}
\newcommand\scalar[1]{\left(#1\right)}
\newcommand\Scalar[1]{\scalar{#1}_{[0,1]}}
\newcommand\Span{\mathop{\rm span}}
\newcommand\supp{\mathop{\rm supp}}
\newcommand\ugauss[1]{\left\lfloor#1\right\rfloor}
\newcommand\with{\, : \,}
\newcommand\Null{{\bf 0}}
\newcommand\bA{{\bf A}}
\newcommand\bB{{\bf B}}
\newcommand\bR{{\bf R}}
\newcommand\bD{{\bf D}}
\newcommand\bE{{\bf E}}
\newcommand\bF{{\bf F}}
\newcommand\bH{{\bf H}}
\newcommand\bU{{\bf U}}
 
\newcommand\cH{{\cal H}}
\newcommand\sinc{{\rm sinc}}
\def\enorm#1{| \! | \! | #1 | \! | \! |}

\newcommand{\dm}{\frac{d-1}{d}}

\let\bm\bf
\newcommand{\balpha}{{\mbox{\boldmath$\alpha$}}}
\newcommand{\bbeta}{{\mbox{\boldmath$\beta$}}}
\newcommand{\bal}{{\mbox{\boldmath$\alpha$}}}
\newcommand{\bbi}{{\bm i}}

\newcommand{\dI}{\Delta}
 
\maketitle
\date{}

\begin{abstract}

 This paper considers the   problem of optimal recovery of an element $u$ of a Hilbert space $\cH$ from measurements of the form $\ell_j(u)$, $j=1,\dots,m$, where the $\ell_j$ are known linear functionals on $\cH$.   Problems of this type are well studied  \cite{MRW} and usually are carried out under an assumption that $u$ belongs to a prescribed  model class, typically a known
 compact subset of $\cH$.   Motivated by reduced modeling for solving parametric partial differential equations, this paper considers another setting where  the additional information about $u$ is in the form of how well $u$ can be approximated by a certain known subspace $V_n$ of $\cH$ of dimension $n$,  or more generally,  in the form of    how well  $u$  can be approximated by each of a sequence of nested subspaces $V_0\subset V_1\cdots\subset V_n$ with each $V_k$ of dimension $k$.   	A  recovery algorithm for the one-space formulation  was proposed in \cite{MPPY}.  Their algorithm is proven,  in the present paper,  to be optimal.  It is also shown how
 the recovery problem for the one-space problem,  has a simple formulation, if  certain favorable bases are chosen to represent $V_n$ and the measurements.  The major contribution of the present  paper is to analyze the multi-space case.  It is shown that, in this multi-space case,  the set of all $u$ that satisfy the given information can be described as the intersection of a family of known ellipsoids in $\cH$.  It  follows    that  a near optimal  
 recovery  
 algorithm in the multi-space problem is provided   by 
 identifying
 any point in  this intersection.   It is easy to see that  the accuracy of recovery of $u$ in the multi-space setting can be much 
  better 
 than in the one-space problems.   Two iterative  algorithms based on alternating projections are proposed for recovery in the multi-space problem and one of them is analyzed in detail.   This analysis includes an a posteriori estimate for the performance of the iterates.  These a 
 posteriori  
 estimates   can serve both as a stopping criteria in the algorithm and also
 as a method to derive convergence rates.
 Since the limit of the algorithm is a point in the intersection of the aforementioned ellipsoids, it provides a near optimal recovery for $u$.       
 \end{abstract}



\section{Introduction}
\label{sec:intro}
 
\subsection{Background and motivation}
The emergence of computational and experimental engineering has led to a spectrum of new mathematical questions on how
 to best merge
  {\em data driven} and {\em model based} approaches.  The development of corresponding {\em data-assimilation} methodologies has been
  originally driven mainly
 by meteorological research (see e.g. \cite{Daley,Lewis}) but has meanwhile entered numerous areas in science and engineering bringing, in particular,
 the role of reduced order modeling into the focus of attention \cite{Bash}.
 
 The present paper addresses some principal mathematical aspects that  arise when trying to numerically capture 
a function $u$ which is a state of a physical process with a known law, however with unknown parameters.
We are given measurements of this state and the question is how to best merge these measurements 
with the model information to come up with a good approximation to $u$. 

A typical setting of this type occurs when  {\em all states}  of the physical process are  described by a specific  {\em parametric family of PDEs} which is known to us, in a form
\begin{eqnarray}
 \nonumber
\cP(u,\mu)=0,
\end{eqnarray}
where $\mu$ is a vector of parameters ranging in a finite or infinite dimensional set $\cP$. Instead of
knowing the exact value of $\mu$ which would allow us to compute the state $u=u(\mu)$
by solving the equation, we observe one of these states
through some collection of measurements and we want to use these measurements, together with the known parametric PDE, to numerically capture the state, or perhaps even more ambitiously to capture the parameters. Since the solution manifold 
\begin{eqnarray}
 \nonumber
 \cM:=\{u(\mu)\; : \; \mu\inÊ\cP\},
 \end{eqnarray}  
to a parametric PDE is generally quite complicated, it is usually seen through a sequence of
nested finite dimensional spaces 
\begin{eqnarray}
 \nonumber
V_0\subset V_1\subset \cdots \subset V_n,\quad  \dim(V_j)=j,
\end{eqnarray}
such that each $V_j$ approximates $\cM$ to a known tolerance $\e_j$. Construction
 of such spaces is sometimes referred to as {\em model reduction}.
 Various algorithms for generating such spaces, together with error bounds $\e_j$, have been derived and analyzed. 
One of the most prominent of these is the {\em reduced basis method} where the spaces
 are generated through particular solution instances $u(\mu^i)$ picked from $\cM$, see \cite{BMPPT,BCDDPW,DPW2,PW}. 
 Other algorithms with known error bounds are based on polynomial approximations in the parametric variable, see \cite{CDS1,CCDS}.

 Thus, the information that the state $u$ we wish to approximate is on the manifold is replaced by the information of how well $u$ can be  approximated by  the   spaces $V_j$.   Of course, this is not enough information to pin down $u$ since we do not know where $u$ is on the manifold, or in the new formulation, which particular element of $V_j$ provides a good approximation to $u$.   However, additional information about $u$ is given by physical measurements  which  hopefully are enough to approximately locate $u$.    This type of recovery problem was formulated and analyzed in \cite{MPPY}
 using an infinite dimensional Hilbert space setting which allows one to properly exploit
 the nature of the continuous background model when assimilating observations. This is also the setting adopted in the present paper. 
 
 The achievements of the present paper are two-fold.   First,  we establish that the algorithm proposed in \cite{MPPY} for estimating a state from a given set of observations and the knowledge of its approximability from a space $V_n$ is {\em best possible in the sense of optimal recovery}.   
 Second, and more importantly, we demonstrate the potential gain in accuracy for state recovery  
 when combining the approximability by {\em each} of the
 subspaces $V_j$ in the given hierarchy. We refer to this  
as the {\em  multi-space setting} which will be seen to better exploit the information given by  reduced bases or polynomial constructions.   We give algorithms and performance bounds for these recovery algorithms in the multi-space setting when the observations are fixed and given to us.
These algorithms are online implementable, similar to the ones discussed in \cite{MPPY}.
Let us mention that one emphasis in \cite{MPPY} is on the selection of the measurement functionals in order to 
optimize the recovery process, while in the present paper we consider such functionals as given and focus on
optimal recovery as explained above.

\subsection{Conceptual preview}
 
We study   the above problems in the general framework  of {\it optimal recovery }  in a Hilbert space  $\cH$ with  inner product $\langle\cdot,\cdot\rangle$ and norm $\|\cdot \|$.
 Under this setting, we are wanting to recover a function $u\in \cH$ from  its measurements $\ell_i(u)=\langle u,\omega_i\rangle$, where the $\omega_i$ are known elements of $\cH$,  $i=1,\dots,m$. If we denote by $W$
the space spanned by the $\omega_i$, $i=1,\dots,m$, then, the measurements determine 
{$w=P_W u$}
 where throughout this paper
$P_X$ denotes the orthogonal projection onto $X$ for any closed subspace $X\subset \cH$.  In going further, we   think of measurements as simply providing  the knowledge of this projection.  In particular, we assume  that  the $\omega_j$'s are linearly independent i.e., $\dim W=m$.    Therefore, our problem is to find an approximation $\hat u(w)$  to $u$ from the information $w\in W$. This is equivalent to constructing a mapping
 $A:W\to \cH$ and setting $\hat u(w)=A(w)=A(P_Wu)$.
 
  All elements of the orthogonal complement $W^\perp$ of $W$ have zero measurements. A first observation is that if all the information we have about $u$ is that 
  {$P_W u= w$}, 
   then we cannot
 recover $u$  to any guaranteed accuracy.  Indeed,  if $u_0$ satisfies the measurements then   $u$ could be any of the   functions $u_0+\eta$, with $\eta\in W^\perp$, and each of these functions would be assigned the same approximation  $\hat u=\hat u(w)$. Therefore, we need additional information about $u$ to have a meaningful problem.   A typical assumption is that $u$ is in some known compact set $\cS\subset \cH$.    The recovery problem  in this case is  known as {\it optimal recovery}.  A classical setting is that $\cH$ is the space $L_2$ and $\cS$ is a finite ball in a Sobolev or Besov space,
 see e.g. \cite{BB, MR, MRW}.

   In contrast to the case  where $\cS$ is a known Sobolev or Besov ball, our interest is in the setting where $\cS$ is the solution manifold $\cM$ of a parametric PDE.    As noted above, the  typical way of resolving $\cM$ is through a finite sequence of spaces $\{V_0,\dots,V_n\}$ with $V_k$ of dimension $k$ where the spaces are known to approximate $\cM$ to some known accuracy.
 This leads us to the following  two settings:
  \vskip .1in
  \noindent
  {\bf The one-space problem:}  We assume that  all what we know about $\cM$ is that there is a  space $  V_n$  of dimension $n$ which is an approximation  to  $\cM$ with accuracy $\e_n$.   Accordingly, we define 
  \be
  \label{defK1}
 \cK:=\cK^{\one}:=\{u\in\cH: 
 { \dist (u,V_n)}
 \le \e_n\},
  \ee
   and consider $u\in\cK$ to be the  only information we have about $\cM$.  In this case, the information  \eref{defK1} is the additional knowledge we have about $u$.  We want to combine this knowledge with our measurements  
   {$P_Wu$}
    to construct a good approximation $\hat u$  to $u$.   So in this case, the spaces $V_n$ and $W$ are known and fixed.
\vskip .1in
  \noindent  
  {\bf The multi-space problem:}  We assume that what we know about $\cM$ is that there is a sequence of spaces $V_0\subset    V_1\subset \cdots\subset V_n$ such that each $V_k$ has dimension $k$ and approximates $\cM$ with accuracy $\e_k$,
  where $\e_0\ge \e_1\ge\cdots\e_n>0$.   This leads us to define  %
  \be
  \label{defKm}
\cK:=\cK^{\mult}:=\bigcap_{j=0}^n \cK^j,
 \ee
 where
  \begin{eqnarray}
 \nonumber
 \cK^j:=\{u\in \cH:  
 {\dist (u,V_j)}
 \le \e_j\}, \quad j=0,\dots,n.
 \label{cKj}
 \end{eqnarray}
 In this case, the information  $u\in \cK$  is the additional knowledge we have about $u$.  We want to combine this knowledge with our measurements to construct a good approximation $\hat u$  to $u$.   As already noted, the multi-space problem is typical when applying reduced bases 
 or polynomial methods to parametric PDEs.

 \subsection{Performance criteria}
  
 This paper is concerned with approximating a function $u\in\cH$ from the information that $u\in \cK$ and $P_Wu=w$
in the two above settings.  Note that in both settings, the set $\cK$ is not compact.
The additional information provided by the measurements gives that 
 $u$ is in the class
 \begin{eqnarray}
 \nonumber
    \cK_w:=\{u\in\cK: P_Wu=w\}.
\end{eqnarray}   
This set is the intersection of $\cK$ with the  affine space
 \begin{eqnarray}
 \nonumber
 \cH_w:=\{u\in\cH: P_Wu=w\}=w+W^\perp.
 \end{eqnarray}
Note that $\cK_w$ may be an empty set for certain $w\in W$.

 Recall that an algorithm is a mapping $A: W\to \cH$ which assigns to any $w\in W$
  the approximation $  \hat u(w)=A(P_wu)$.    In designing an algorithm, we are given the information
  of the spaces $(V_k)_{k=0,\dots,n}$ and the error bounds $(\e_k)_{k=0,\dots,n}$.
  There are several ways in which we can measure the performance of  an algorithm.    Consider first the one-space problem.     A first way of measuring the performance of an algorithm is to ask for 
 an estimate of the form
 \be
  \label{instopt}
  \|u-A(P_Wu)\|\le C_A(w)\dist(u,V_n),\quad u\in\cK_w.
  \ee
 The best algorithm $A$, for a given fixed value of $w$, would give the smallest constant $C_A(w)$ and the algorithm which gives this smallest constant is said to be {\it instance optimal} with constant $C_A(w)$.   In this case, the performance bound given by the right side of \eref{instopt}
  depends not only on $w$ but on the particular $u$ from $\cK_w$.

 The estimate \eref{instopt} also gives a performance bound for the entire class $\cK_w$ in the form
   \begin{eqnarray}
 \nonumber
  \sup_{u\in\cK_w} \|u-A(P_Wu)\|\le C_A(w)\e_n.
  \end{eqnarray}
 This leads us to the notion of performance of a recovery algorithm $A$ on any set $\cS\subset \cH$ which is  defined by  
 \begin{eqnarray}
 \nonumber
E_A(\cS):= \sup_{u\in \cS} \|u-A(P_Wu)\|.
  \end{eqnarray}
  The {\it class optimal performance} on the set  $\cS$ is given by
  \be
  \label{optalg}
  E(\cS):=\inf_{A}E_A(\cS),
  \ee
  where the infimum is taken over all possible algorithms, i.e., all maps $A:W\to \cH$.
  In particular,   class optimal performance is defined for  both the single space or multi-space settings and for both
  the sets $\cK_w$ for each of individual $w$ which gives the measure $E(\cK_w)$ or  
  the entire class $\cK$ which gives the performance $E(\cK)$. The latter notion is the most meaningful when  in applications it is not known
  which measurements $w\in W$ will appear or will be available.
     
 The present paper studies each of the above problems with the goal of determining the best algorithms. 
 For this purpose, we 
 introduce for any closed subspaces $V$ and $W$ of $\cH$ the quantity
 \be
\label{mu}
 \mu(V,W):= \sup_{\eta\in W^\perp} \frac{\|\eta\|}{\|\eta-P_V\eta\|}=\sup_{\eta\in W^\perp} \frac{\|\eta\|}{\|P_{V^\perp} \eta\|}.
\ee
A simple calculation shows that $\mu(V,W)=\beta(V,W)^{-1}$ where
 \begin{eqnarray}
\nonumber
\beta(V,W):=\inf_{v\in V}\frac{\|P_Wv\|}{\|v\|} =\inf_{v\in V}\sup_{w\in W}\frac{ \langle v,w\rangle}{\|v\|\|w\|} .
\end{eqnarray}
Note that in the case where $V=\{0\}$ we have $\mu(V,W)=1$.
  
In \S \ref{sec:one1}  of the paper, we analyze the one space problem, that is, $\cK=\cK^{\rm one}$.  
The inf-sup constant $\beta$ was used in \cite{MPPY}, for the study of this  problem,
where the authors proposed an algorithm, in the form of a certain linear mapping $A^*:w\to A^*(w)$, then analyze its performance.    
While the approach in \cite{MPPY}
  is based on variational arguments, ours is quite different and geometric in nature.   Our first goal is to establish that the algorithm proposed in
  \cite{MPPY} is both instance optimal and class optimal. We show that  for any function $u\in\cH$   %
  \be
  \label{errorbound1}
  \|u-A^*(P_Wu)\|\le \mu(V_n,W)\dist(u,V_n) .
  \ee
  Notice that if $\beta(V_n,W)=0$,  the above estimate would give no bound on approximation as is to be
  expected since $V_n$ would contain elements of $W^\perp$ and these cannot be distinguished by the measurements.  This would always be the case if $n>m$ and so in going further we always work under the assumption that $n\le m$.  
 
  Let us note that this is a modest improvement on the estimate in \cite{MPPY} 
  which has the constant $\mu(V_n,W)+1$ rather than $\mu(V_n,W)$ on the right side of \eref{errorbound1}.  
  More importantly, we show that the estimate \eref{errorbound1} is best possible in the sense that the constant $\mu(V_n,W)$ cannot be replaced by a smaller
   constant.   Another important remark, observed in \cite{MPPY},  is that in \eref{errorbound1},  $\dist(u, V_n)$ can be replaced by the smaller quantity $\dist(u,V_n\oplus(W\cap V_n^\perp))$.   We establish, with our approach, the estimate  
  \be
  \label{errorbound12}
  \|u-A^*(P_Wu)\|\le \mu(V_n,W)\dist(u,V_n\oplus(W\cap V_n^\perp)) ,
  \ee
  which improves the constant given in \cite{MPPY}.
We again show that $\mu(V_n,W)$ is the best constant in estimates of this form.
  
In view of  \eref{errorbound1}, the algorithm $A^*$ provides the class estimate 
  \be
  \label{errorbound13}
 E_{A^*}(\cK)\le    \mu(V_n,W)\e_n.
  \ee
  We again show that this  algorithm is {\it class optimal} in the sense that for the single space problem
 \begin{eqnarray}
\nonumber
 E (\cK)=    \mu(V_n,W)\e_n.
  \end{eqnarray}
Our analysis is based on proving lower bounds which show that the upper 
{estimates  \eref{errorbound12}}
 and \eref{errorbound13} cannot be improved.
These lower bounds apply to both linear and nonlinear algorithms, that is, 
 \eref{errorbound12} and \eref{errorbound13} cannot be improved also using
 nonlinear mappings.

Another goal of our analysis of the one-space problem is to simplify the description of the optimal solution 
through the choice of, what we call, {\em favorable bases} for
the spaces $V_n$ and $W$.   These favorable bases are then used in our analysis of the multi-space problem
which is the object of \S \ref{sec:multiple}. One  possible way of proceeding, in the multi-space case,  is to examine
  the right side of \eref{errorbound13} for each of the spaces $(V_k)_{k=0,\dots,n}$, and choose the one which gives the minimum  value.   This would produce an algorithm $A$ with the error bound
  \be 
  \label{olderror}
E_A(\cK)\le   \min_{0\le k\le n} \mu(V_k,W)\e_k.
  \ee
   Notice that the $\e_k$ are decreasing but
  the $\mu(V_k,W)$ are increasing as $k$ gets larger.   So these two quantities are working against one another and the minimum may be assumed for  an intermediate value of $k$.

  It turns out that the algorithm giving the bound \eref{olderror}  may be far from optimal and our main achievements 
  in \S \ref{sec:multiple} are to produce both algorithms and  a priori  performance bounds
  which in general  are better than that of \eref{olderror}.  We show how the multi-space problem is connected to finding
  a point in the intersection of a family of ellipsoids in $\cH$ and propose an algorithm  based on this intersection  property.
  Then,  we give {\em a priori bounds} on the performance of our numerical algorithm,
  which are shown to be, in general, better than  \eref{olderror}.

\section{The one-space problem}
\label{sec:one1}

\subsection{Preliminary remarks}\label{sec:remarks}
 \label{sec:prel}

We begin with some general remarks which can be applied to our specific problem.   If $\cS\subset \cH$ is a bounded set and we wish to simultaneously 
approximate  
all of the elements in $\cS$,
 then the best approximation is described by  the center of the {\it Chebyshev  ball} of $\cS$, which is defined as the smallest  closed ball that contains
$\cS$.  To describe this ball, we first define    the {\em Chebyshev  radius}
 \begin{eqnarray}
\nonumber
\rad(\cS) :=  \inf \{r: \cS\subset B(v,r)\  {\rm for \ some} \ v\in\cH\}.
\end{eqnarray}
The following well known lemma says that the Chebyshev  ball  exists and is unique. 

\begin{lemma}
\label{lemma:radassumed}
If  $\cS$ is any bounded set in $\cH$  with $R:=\rad(\cS)$,  then there exists a unique $v^*\in \cH$ such that
 \be
\label{assumed}
\cS\subset B(v^*, R).
\ee
\end{lemma}

\noindent
{\bf Proof:} 
 For any $v\in\cH$, we define
\begin{eqnarray}
\nonumber
R_\cS(v):= \inf\,\{r: \cS\subset B(v,r)\},
\end{eqnarray}
which is  a well-defined function from $\cH$ to $\R$.
It follows from triangle inequality that  $R_\cS : \cH \to \R$ is continuous. It is also easily seen that
\begin{eqnarray}
\nonumber
\cS\subset B(v,R_\cS(v)).
\end{eqnarray}
By definition, $\rad(\cS)= \inf_{v\in\cH}
R_{\cS}(v)$.    Now,  consider any infimizing sequence $(v_j)_{j\in\N}$, i.e.,
\begin{eqnarray}
\nonumber
\lim_{j\to\infty} R_{\cS}(v_j)= \rad(\cS).
\end{eqnarray}
We claim that $(v_j)_{j\in\N}$ is a Cauchy sequence. To see this,  define   $r_j:= R_{\cS}(v_j)$.  
For any fixed $j$ and $k$ and any $z \in \cS$  we  define 
  $d_j:= v_j-z$ and $d_k:=v_k-z$.    Then, $\|d_j\|\le r_j$, and $\|d_k\|\le r_k$.   Therefore,
\begin{eqnarray*}
\|v_j - v_{k}\|^2 &=&\|d_j -d_{k}\|^2 = \langle d_j -d_{k},d_j -d_{k}\rangle\\
&=& 2\langle d_j,d_j\rangle + 2 \langle d_{k},d_{k}\rangle -\langle d_j + d_k,d_j+ d_k\rangle\\
&=& 2 \|d_j\|^2 + 2 \|d_k\|^2 - 4\Big\|\frac 12 (d_j + d_k)\Big\|^2\\
&\leq& 2 r_j^2 + 2 r_k^2 -4\Big\|\frac 12 (v_j + v_k)-z\Big\|^2.
\end{eqnarray*}
Since $z\in \cS$ is arbitrary we get
$$
\|v_j - v_{k}\|^2 \leq 2 r_j^2 + 2 r_k^2 -4 \[ R_{\cS}\( \frac 12 (v_j+v_k)\) \]^2 \leq  2 r_j^2 + 2 r_k^2 -4\rad(\cS)^2.
$$
Since $r_j, r_k \to\rad(\cS)$,  this shows that $(v_j)_{j\in\N}$ is a Cauchy sequence and has a limit $v^*$,  which
by  the continuity of $v\mapsto R_{\cS}(v)$ satisfies $R_{\cS}(v^*)=\rad(\cS)$. The uniqueness of $v^*$
also follows from the above inequality by contradiction. By using the continuity of $v\mapsto R_{\cS}(v)$ one easily shows that \eref{assumed} holds.
\hfill $\Box$
\nl

We sometimes say that $v^*$ in the above lemma is the {\it center} of $\cS$. 
For any bounded set $\cS$, the diameter of $\cS$  is related to its Chebyshev radius $\rad(\cS)$ by the inequalities
\begin{eqnarray}
\nonumber
\rad(\cS)\le  \diam(\cS)\le 2 \rad(\cS). 
\end{eqnarray}
For general sets $\cS$  these inequalities cannot be improved.    However,  we have the following   remark.
\begin{remark}
\label{rem:diam2r}
Let $\cS$ be symmetric about a point $z$, i.e. whenever $v\in \cS$, then $2z-v\in \cS$.  Then,
the Chebyshev radius of $\cS$ equals half its diameter, that is, $\diam(\cS)=2\rad(\cS)$ and its center is $z$. 
\end{remark}

\begin{remark}
\label{rem:opt}
In the particular setting of this paper, for any given $w\in W$ such that $\cK_w$ is non-empty,
the optimal recovery $u^*(w)$ over the class $\cK_w$ is obviously 
given by the center of $\cK_w$, and the class optimal performance is  given by
\begin{eqnarray}
\nonumber
E (\cK_w)=\rad(\cK_w).
\end{eqnarray}
\end{remark}
 \begin{remark}
\label{rem:center}
For a bounded, closed, convex set $\cS\subset \cH$ (which is always the case in this paper) its center  $u$ is in $\cS$. In fact, if  this was not true,
by translating $\cS$, we can assume $u=0$. Let $s_0=\argmin_{s\in \cS}\|s\|$. By convexity $s_0$ exists, $s_0\neq 0$, and $\<s,s_0\>\geq \<s_0,s_0\>$,
$s\in \cS$. Thus 
$$
\sup_{s\in \cS}\|s-s_0\|^2=\sup_{s\in \cS} ( \<s,s\>-2\<s,s_0\>+\<s_0,s_0\>)\leq \sup_{s\in \cS} \|s\|^2-\|s_0\|^2
$$
which contradicts the assumption that $0$ is the center of $\cS$.
\end{remark}

\subsection{Optimal bounds for the one-space problem}
\label{sec:one}

We next consider the case where the set 
$\cK=\cK^{\one}$ is given by \iref{defK1},
where $V_n$ is a fixed and known $n$ dimensional space. 
 In this section, we derive the algorithm proposed in \cite{MPPY}, however from a different
point of view emphasizing more the optimal recovery and geometric aspects of the problem.   This
allows us to improve on their estimates some but, more importantly, it is also useful when treating the multi-space problem. 
 \nl
 
{\em In the event that  $\beta(V_n,W)=0$, the space $V_n$ contains elements from 
$W^\perp$ which implies that if $w\in W$ is such that $\cK_w$ is non-empty, then
$\cK_w$ is unbounded, or equivalently $\rad(\cK_w)$ is infinite, which means
that we cannot hope for any guaranteed performance over $\cK_w$. This is the case
in particular when $n>m$. For this reason, in the rest of the paper, we always assume
that $\beta(V_n,W)>0$, which means in particular that $n\le m$.}
\nl

Let $w$ be any element from $W$.
We claim that the map
\begin{eqnarray}
\nonumber
u\mapsto  \|u-P_{V_n}u\|= \|P_{V_n^\perp}u\|,
\end{eqnarray}
admits a unique minimizer over the affine space $\cH_w$. To see this, we let  $u_0$ be any element from  $\cH_w$.  It follows  that
every $u\in \cH_w$ can be written as
$u=u_0+\eta$ for some $\eta \in W^\perp$.
Minimizing $\|P_{V_n^\perp}u\|$ over $\cH_w$ therefore amounts to minimizing the function 
\begin{eqnarray}
\nonumber
\eta \mapsto  f (\eta):=\|P_{V_n^\perp}u_0+P_{V_n^\perp}\eta\|^2,
\end{eqnarray}
over $W^\perp$. We may write
\begin{eqnarray}
\nonumber
f(\eta):=g(\eta)+\|P_{V_n^\perp}\eta\|^2,
\end{eqnarray}
where $g$ is an affine function. Since we have assumed that  $\beta(V_n,W)>0$, the inequalities
\begin{eqnarray}
\nonumber
\beta(V_n,W) \|\eta\| \leq \|P_{V_n^\perp}\eta\|Ê\leq \|\eta\|,\quad \eta \in W^\perp.
\end{eqnarray}
show that $\eta\mapsto  \|P_{V_n^\perp}\eta\|$ is an equivalent norm over $W^\perp$. Therefore
$\eta \mapsto f(\eta)$ is strongly convex over $W^\perp$ and therefore admits a unique minimizer 
\begin{eqnarray}
\nonumber
\eta^* :=\argmin_{\eta\in W^\perp} f(\eta).
\end{eqnarray}
It follows that $u^*=u_0+\eta^*$ satisfies
 \begin{eqnarray}
\nonumber
u^*=u^* (w) :=\argmin_{u\in\cH_w} \|u-P_{V_n}u\|
\end{eqnarray}
and that this minimizer is unique.

\begin{remark}
\label{rem:nonempty}
If $w$ is such that $\cK_w$ is non-empty, there exists a $u\in \cH_w$ such that $\|u-P_{V_n}u\|\leq \e_n$.
Therefore $\|u^*-P_{V_n}u^*\|\leq \e_n$, that is, $u^*\in\cK_w$.  In particular, $u^*$ minimizes $\|u-P_{V_n}u\|$ over all $u\in \cK_w$.
\end{remark}
We next define 
\begin{eqnarray}
\nonumber
v^*:=v^*(w):=P_{V_n}u^*.
\end{eqnarray} 
From the definition of $u^*$, it follows that the pair $(u^*,v^*)$ is characterized by the minimization property
\be
\|u^*-v^*\|=\min_{u\in \cH_w, \,v\in V_n} \|u-v\|,
\label{doublemin}
\ee
As the following remark shows, $u^*-v^*$ has a certain double orthogonality property.

 \begin{remark}
 \label{rem:orthogonal}  The element  $u^*-v^*$ is orthogonal to both spaces $V_n$ and $W^\perp$.
The orthogonality to $V_n$ follows from the fact that $v^*=P_{V_n}u^*$. On the other hand, for any
$\eta\in W^\perp$ and $\alpha\in \R$, we have
\begin{eqnarray}
\nonumber
\|u^*-v^*\|^2\leq \|u-P_{V_n}u\|^2,\quad u:=u^*+\alpha \eta,
\end{eqnarray}
and thus
\begin{eqnarray}
\nonumber
\|u^*-v^*\|^2\le \|u^*-v^*+\alpha(\eta -P_{V_n}\eta) \|^2= \|u^*-v^*\|^2+2\alpha \langle u^*-v^*,\eta\rangle +\alpha^2\|\eta -P_{V_n}\eta\|^2.
\end{eqnarray}
This shows that $u^*-v^*$ is orthogonal to $W^\perp$.
\end{remark}

 \begin{remark}
 \label{rem:uniquepair} Conversely, if $u\in \cH_w$ and $v\in V_n$ are such that  $u-v$ is orthogonal to both spaces $V_n$ and $W^\perp$,
 then $u=u^*$ and $v=v^*$.  Indeed,  from  this orthogonality
 \begin{eqnarray}
\nonumber
\|u^*-v^*\|^2 = \|u-v\|^2+\|u^*-v^*-(u-v) \|^2.
\end{eqnarray}
 This gives that $u,v$ is also a minimizing pair and from uniqueness of the minimizing pair $u=u^*$ and $v=v^*$.
\end{remark}

The next theorem describes the smallest ball that contains $\cK_w$, i.e., the Chebyshev ball for this set,
and shows that the center of this ball is $u^*(w)$.

\begin{theorem} 
\label{theo:BA} Let $W$ and $V_n$ be such that $\beta(V_n,W)>0$. 

\noindent
{\rm (i)}
For any $w\in W$ such that $\cK_w$ is non-empty, the Chebyshev ball for $\cK_w$ is the ball centered at $u^*(w)$ of radius 
\be
\label{Rstar}
R^*=R^*(w):=\mu(V_n,W) (\e _n^2-\|u^*(w)-v^*(w)\|^2)^{1/2}. 
\ee

\noindent
{\rm (ii)} The optimal algorithm in the sense of \eref{optalg} for recovering  $\cK_w$ from the measurement $w$ is given by the mapping $A^*: w\mapsto u^*(w)$ and gives the performance bound
 \be
\label{perfba}
E_{A^*}(\cK_w)= E(\cK_w)=\mu(V_n,W)(\e_n^2-\|u^*(w)-v^*(w)\|^2)^{1/2}.
\ee
 {\rm (iii)} The optimal algorithm in the sense of \eref{optalg} for recovering  $\cK$ is given by the mapping $A^*: w\mapsto u^*(w)$ and gives the performance bound
 \be
\label{perfba1}
E_{A^*}(\cK)=E(\cK)=\mu(V_n,W)\e_n.
\ee
 \end{theorem}

\noindent
{\bf Proof:}   In order for $\cK_w$ to be nonempty, we need that $\|u^*-v^*\|\le \e_n$.    Any $u\in\cH_w$  can be written as  $u=u^*+\eta$ where $\eta\in W^\perp$.    Therefore,
  \begin{eqnarray}
\nonumber
 u-P_{V_n}u=
u^*-v^*+\eta-P_{V_n}\eta.
\end{eqnarray}
Because of the orthogonality in Remark \ref{rem:orthogonal}, we have
\be
\label{min3}
  \|u-
  {P_{V_n}u}
  \|^2=\|u^*-v^*\|^2+\|\eta-
  {P_{V_n}\eta}
  \|^2.
\ee
Thus a necessary and sufficient condition for $u$ to be in $\cK_w$ is that
 \begin{eqnarray}
\nonumber
 \|P_{V_n^\perp}\eta\|^2= \|\eta-
 {P_{V_n}\eta}
 \|^2\le \e_n^2-\|u^*-v^*\|^2.
\end{eqnarray}
From the definition of $\mu(V_n,W)$, this means that any $u\in\cK_w$ is contained in the ball $B(u^*(w),R^*(w))$.     Now, if $\eta$ is any element   in $W^\perp$  with norm $R^*(w)$ which achieves
the maximum in the definition of $\mu(V_n,W)$, then $u^*\pm \eta$ is in $\cK_w$ and since $\|\eta\| =R^*(w)$ we see that the diameter of $\cK_w$ is at least as large as $2R^*(w)$.  Since $\cK_w$ is the translation of a symmetric set, we thus obtain (i) from Remark \ref{rem:diam2r}.   The claim (ii) about $A^*$ being the  optimal algorithm follows from Remark \ref{rem:opt}.  Finally, the 
performance bound \iref{perfba1} in the claim (iii) holds because the maximum of $R^*(w)$ is achieved when $w=0$.
\hfill $\Box$

\begin{remark}
The optimal mapping $w\mapsto A^*(w)=u^*(w)$ is independent of $\e_n$
and the knowledge of $\e_n$ is not needed in order to compute $A^*(w)$.
\end{remark}

\begin{remark} 
\label{rem:character}
Since $\cK_w$ is the intersection of the cylinder $\cK$ with the affine space $\cH_w$, it has
the shape of an ellipsoid. The above analysis describes this ellipsoid as follows:
a point $u^*+\eta$ is in $\cK_w$ if and only if $\|P_{V_n^\perp}\eta\|^2\le  \e_n^2-\|u^*-v^*\|^2$. 
In the following section, we give a parametric description of this ellipsoid using certain coordinate systems,
see Lemma \ref{lem:ellipsoids}.
\end{remark}

\begin{remark}
\label{rem:MP}  The elements $u^*$ and $v^*$ were introduced in \cite{MPPY} and used 
to define the algorithm $A^*$ given in the above theorem.  The analysis from \cite{MPPY} establishes the error bound
 \begin{eqnarray}
\nonumber
\|u-u^*(w)\|\le (\mu(V_n,W)+1)\dist(u,V_n\oplus (V_n^\perp\cap W)).
\end{eqnarray}
A sharper form of this inequality can be derived from our results.   Namely,   if $u$ is any element in $\cH$ then we can define $\e_n:=\|u-P_{V_n}u\|$ and $w:=P_Wu$.  Then,  $u\in \cK_w$,  for this choice of $\e_n$, and so Theorem \ref{theo:BA}  applies
and gives a recovery of $u$ with the bound 
\be
\label{min5}
\|u-u^*(w)\|\le \mu(V_n,W)(\e _n^2-\|u^*-v^*\|^2)^{1/2}
=  \mu(V_n,W)\|u-P_{V_n} u-(u^*-v^*)\|,
\ee
where the second equality follows from \iref{min3}. We have noticed in Remark \ref{rem:orthogonal} that $u^*-v^*\in V_n^\perp\cap W$,
and on the other hand we have that 
{$u-(u^*-v^*) \in V_n+W^\perp$}, 
which shows that 
\begin{eqnarray}
\nonumber
u^*-v^*=P_{V_n^\perp\cap W}u.
\end{eqnarray}
Therefore 
\begin{eqnarray}
\nonumber
P_{V_n} u+u^*-v^*=P_{V_n\oplus (V_n^\perp\cap W)} u,
\end{eqnarray}
and \iref{min5} gives
\begin{eqnarray}
\nonumber
\|u-u^*(w)\|\le \mu(V_n,W)\dist(u,V_n\oplus (V_n^\perp\cap W)).
\end{eqnarray}
\end{remark}

\begin{remark}
\label{rem1}
 Let us observe that given a space $V_n$ with   $n<m$ we have $(W\cap V_n^\perp)\neq \{0\}$,  thus the space  $\bar V_n:=
V_n\oplus (W\cap V_n^\perp)$ is strictly larger than $V_n$. However $\mu(\bar V_n,W)= \mu( V_n,W) $ because for any $\eta\in W^\perp$, the projection of $\eta$ onto $W\cap V^\perp_n$ is zero. In other words we can enlarge $V_n$ preserving the estimate \iref{perfba1} for class optimality performance as long as we add parts of $W$ that are
orthogonal to $V_n$.
\end{remark}

  \subsection{The numerical implementation of the optimal algorithm}
 \label{ss:numericalone-space} 
  Let us next discuss the numerical implementation of the optimal algorithm for the {one-space}
  problem.   Let $\omega_1,\dots,\omega_m$ be any orthonormal basis for $W$.  For  theoretical reasons only, we complete it to an orthonormal basis for $\cH$.  So   $\{\omega_i\}_{i>m}$ is a complete orthonormal system for $W^\perp$.  We can write down explicit formulas for $u^*$ and $v^*$.   Indeed, any $u\in   \cH_w$ can be written
   \begin{eqnarray}
\nonumber
u= \sum_{i=1}^m w_i\omega_i+\sum_{i=m+1}^\infty x_i\omega_i,
 \end{eqnarray}
 where $w_i:=\langle w,\omega_i\rangle$,  and $(x_i)_{i>m}$ is any $\ell_2$ sequence.  So, for any $v\in V_n$ and $u\in  \cH_w$, we have
\begin{eqnarray}
\nonumber
\|u-v\|^2= \sum_{i=1}^m(w_i-v_i)^2+\sum_{i=m+1}^\infty (x_i-v_i)^2,
 \end{eqnarray}
 where $v_i:=\langle v,\omega_i\rangle$.
Thus, for any $v\in V_n$, its best approximation $u(v)$ from $  \cH_w$ is 
\be
u(v):=\sum_{i=1}^mw_i\omega_i+\sum_{i=m+1}^\infty v_i \omega_i,
\label{uv}
\ee
and its error of approximation is 
\begin{eqnarray}
\nonumber
\|v-u(v)\|^2=\sum_{i=1}^m(w_i-v_i)^2.
\end{eqnarray}
In view of \iref{doublemin} we have 
\begin{eqnarray}
\nonumber
v^*=\argmin_{v\in V_n}\|v-u(v)\|^2  =\argmin_{v\in V_n}\sum_{i=1}^m(w_i-v_i)^2=\argmin_{v\in V_n}\|w-P_Wv\|^2.
\end{eqnarray}
For any given orthonormal basis $\{\phi_1,\cdots,\phi_n\}$ for $V_n$, we can find the coordinates of 
$v^*\in V_n$ in this basis by solving the $n\times n$ linear system associated to the above least squares problem.  Once $v^*$ is found,
the optimal recovery $u^*=u^*(w)$ is given, according to \iref{uv}, by 
\begin{eqnarray}
\nonumber
u^*=v^* +\sum_{i=1}^m(w_i-v_i^*)\omega_i,
\end{eqnarray}
where $v_i^*=\<v^*,
{\omega_i}\>$. Note that we may also write
 \be
 \label{lsu}
 u^*=\sum_{i=1}^mw_i\omega_i+\sum_{i=m+1}^\infty \langle v^*,\omega_i\rangle \omega_i = w+P_{W^\perp}v^*.
 \ee

\subsection{Liftings and   favorable bases for $V_n$ and $W$}
\label{ss:favorable}

It turns out that the above optimal algorithm has an even simpler description if we  choose suitable bases for $V_n$ and $W$,
which we call {\em favorable bases}.  These bases will also be important in our analysis of the multi-space problem. To describe this new geometric view, we introduce the  description of algorithms through  liftings and see how the best algorithm of the previous section arises in this context. 

  As noted earlier, any algorithm is  a mapping $A:W\to\cH$ which takes   $w=P_Wu$ into $\hat u(w)=A(w)=A(P_Wu)$.  This image  serves as the approximant of all of the $u\in\cK_w$.   We can write any $u\in \cK_w$ as $u=w+P_{W^\perp} u$.
So the problem is to find an appropriate mapping $F:W\to {W^\perp}$ and take as the approximation
\begin{eqnarray}
\nonumber
\hat u (w):= A(w):=w+F(w).
\end{eqnarray}
 At this stage $F$ can be any  linear or nonlinear mapping from $W$ into $W^\perp$.   We call such mappings $F$ {\it liftings}.

 According to \eref{lsu},  the optimal lifting $F^*$ is defined by
\begin{eqnarray}
\nonumber
F^*(w)= P_{W^\perp}v^*(w)\in P_{W^\perp}V_n,
\end{eqnarray}
which is actually a linear mapping since $v^*$ depends linearly on $w$.   The algorithm $A^*(w)=w+F^*(w)$ was shown in the previous section to be optimal for each class $\cK_w$ as well as for $\cK$. Note that this optimality holds even if we open the competition to nonlinear maps $F$, respectively $A$.

We  next show that $F^*$ has a simple description as a linear mapping by introducing favorable bases.
  We shall make use of the following elementary facts from linear algebra: if $X$ and $Y$ are closed subspaces of  a Hilbert space $\cH$,
then:
\begin{itemize}
\item
We have the equality
\begin{eqnarray}
\nonumber
 \dim({P_X Y})=\dim({P_YX}).
\end{eqnarray}
This can be seen by introducing the cross-Gramian matrix $G=(\<x_i,y_j\>)$, where $(x_i)$ and $(y_j)$ are orthonormal
bases for $X$ and $Y$.  Then $G$ is the  matrix representation of 
 the projection operator $P_X$ from $Y$ onto $X$
with respect to these bases and $G^t$ is the corresponding representation of 
 the projection operator $P_Y$ from $X$ onto $Y$.
Hence,
\begin{eqnarray}
\nonumber\dim(  P_X Y)={\rm rank}(G)={\rm rank}(G^t)=\dim(P_Y X).
\end{eqnarray}
\item
The space $Y$ can be decomposed into a direct orthogonal sum
\be
Y= P_Y X \oplus (Y\cap X^\perp).
\label{sum}
\ee
For this, we need to show that  $Y\cap X^\perp=Z$ where $Z\subset Y$ is the orthogonal complement of {$P_YX$}
 in $Y$. If $y\in Z$, then   $\<y,P_Yx\>=0$ for all $x\in X$.  
Since $\<y,x-P_Yx\>=0$, if follows that $\<y,x\>=0$, for all $x\in X$, and thus $y\in Y\cap X^\perp$. Conversely
if $y\in Y\cap X^\perp$, then for any $x\in X$ 
{$\<y,P_Yx\>=-\<y,x-P_Yx\>=0$,}
 which shows that $y\in Z$.
\end{itemize}

Now to construct the favorable bases we want, we begin with any  orthonormal 
{basis} 
$\{\phi_1,\dots,\phi_n\}$ of  $V_n$ and any orthonormal basis $\{\omega_1,\dots,\omega_m\}$ of $W$.
We consider the $m\times n$ cross-Gramian matrix
\begin{eqnarray}
\nonumber
G:=(\<\omega_i,\phi_j\>),
\end{eqnarray}
which may be viewed as the matrix representation of the projection operator $P_W$ from $V_n$ onto $W$
using these bases since $P_W(\phi_j)=\sum_{i=1}^m\<\omega_i,\phi_j\> \omega_i$. Note that the inf-sup condition  $\beta(V_n,W)>0$ means that
\begin{eqnarray}
\nonumber
\dim({P_W V_n})=n,
\end{eqnarray}
or equivalently, the rank of $G$ is equal to $n$. We perform a singular value decomposition of $G$, which gives
\begin{eqnarray}
\nonumber
G=USV^t
\end{eqnarray}
where $U=(u_{i,j})$ and $V=(v_{i,j})$ are unitary $m\times m$ and $n\times n$ matrices, respectively,
and where $S$ is an $m\times n$ matrix with entries $s_i>0$ on the diagonal $i=j$, $i=1,\dots,n$,  and zero  entries elsewhere.
This allows us to define new  orthonormal bases $\{\phi_1^*,\dots,\phi_n^*\}$ for $V_n$ and
$\{\omega_1^*,\dots,\omega_m^*\}$ for $W$ by
\begin{eqnarray}
\nonumber
\phi_j^*=\sum_{i=1}^n v_{i,j}\phi_i\quad {\rm and} \quad \omega_j^*=\sum_{i=1}^m u_{i,j}\omega_i.
\end{eqnarray}
These new bases are such that
\begin{eqnarray}
\nonumber
P_W(\phi_j^*)=s_j \omega_j^*, \quad j=1,\dots, n,
\end{eqnarray}
and have diagonal
cross-Gramian, namely
\begin{eqnarray}
\nonumber
\<\omega_i^*,\phi_j^*\>= s_j\delta_{i,j}.
\label{cross}
\end{eqnarray}
Therefore 
$\{\omega_1^*,\dots,\omega_n^*\}$ and $\{\omega_{n+1}^*,\dots,\omega_m^*\}$
are orthonormal bases for the
$n$-dimensional space 
{$P_W V_n$}
 and respectively its orthogonal complement in $W$
which is $V_n^\perp\cap W$ according to \iref{sum}. 

By convention, we organize the singular values in decreasing order 
\begin{eqnarray}
\nonumber
0<s_n\leq s_{n-1} \leq \cdots \leq s_1.
\end{eqnarray}
Since $P_W$ is an orthogonal projector, all of them are at most $1$
and in the event where 
\begin{eqnarray}
\nonumber
s_1=s_2=\cdots=s_{p}=1,
\end{eqnarray}
for some $0<p\leq n$, 
then we must have
\begin{eqnarray}
\nonumber
\omega_j^*=\phi_j^*,  \quad j=1,\dots,p.
\end{eqnarray}
This corresponds to the case where $V_n\cap W$ is non-trivial
and $\{\omega_1^*,\dots,\omega_p^*\}$ forms an orthonormal basis for $V_n\cap W$.
We define $p=0$ in the case where $V_n\cap W=\{0\}$. 

We may now give a simple description of
the optimal algorithm $A^*$ and lifting $F^*$, in terms of their action
on the basis elements $\omega_j^*$. For $j=n+1,\dots,m$,
we know that $\omega_j^*\in V_n^\perp\cap W$. From Remark \ref{rem:uniquepair}, it follows that the
optimal pair $(u^*,v^*)$ which solves \iref{doublemin} for $w=\omega_j^*$ is
\begin{eqnarray}
\nonumber
u^*=\omega_j^*\quad {\rm and} \quad v^*=0,
\end{eqnarray}
and therefore
\begin{eqnarray}
\nonumber
A^*(\omega_j^*)=\omega_j^*\quad {\rm and}\quad F^*(\omega_j^*)=0
,\quad j=n+1,\dots,m.
\end{eqnarray}
For $j=1,\dots,n$, we know that $\omega_j^*=P_W(s_j^{-1}\phi_j^*)$. It follows that the
optimal pair $(u^*,v^*)$ which solves \iref{doublemin} for $w=\omega_j^*$ is
\begin{eqnarray}
\nonumber
u^*=v^*=s_j^{-1}\phi_j^*.
\end{eqnarray}
Indeed, this follows from Remark \ref{rem:uniquepair} since this  pair has $u^*-v^*=0$ and hence has the double orthogonality property.
So, in this case,
\begin{eqnarray}
\nonumber
A^*(\omega_j^*)=s_j^{-1}\phi_j^*\quad {\rm and}\quad F^*(\omega_j^*)=s_j^{-1}\phi_j^*-\omega_j^*.
\end{eqnarray}
Note in particular that $F^*(\omega_j^*)=0$ for $j=1,\dots,p$.

 \begin{remark}
  \label{infsuprem}
The favorable bases are useful when computing the inf-sup constant $\beta(V_n,W)$.  Namely,
 for an element $v=\sum_{j=1}^n v_j\phi_j^*\in V_n$ 
  we find that $P_Wv=\sum_{j=1}^n s_jv_j\omega_j^*$
  and so
  \begin{eqnarray}
\nonumber
  \beta(V_n,W)=\min_{v\in V_n} \frac{\|P_Wv\|}{\| v\|}=\min_{v\in V_n} \(\frac{\sum_{j=1}^n s_j^2 v_j^2}{\sum_{j=1}^n  v_j^2}\)^{1/2}=\min_{j=1,\dots,n} s_j=s_n. 
  \end{eqnarray}
  Correspondingly,
    \begin{eqnarray}
\nonumber
   \mu(V_n,W)=s_n^{-1}.
   \label{munew}
   \end{eqnarray}
   Recall that for the trivial space $V_0=\{0\}$, we have $\mu(V_0,W)=1$.
   \end{remark}

For further purposes, we complete the favorable bases into orthonormal bases 
of $\cH$ by constructing particular orthonormal bases for $V_n^\perp$ and $W^\perp$. According to \iref{sum}
we may write these spaces as direct orthogonal sums
\begin{eqnarray}
\nonumber
V_n^\perp=P_{V_n^\perp} (W) \oplus (V_n^\perp \cap W^\perp),
\end{eqnarray}
and
\begin{eqnarray}
\nonumber
W^\perp=P_{W^\perp}(V_n) \oplus (V_n^\perp \cap W^\perp).
\end{eqnarray}
The second space $V_n^\perp \cap W^\perp$ in the above decompositions may be of infinite dimension and we consider 
an arbitrary orthonormal basis $(\psi_i^*)_{i\geq 1}$ for this space. For the first spaces in the above decompositions,
we can build orthonormal bases from the already constructed favorable bases.

For the space $P_{V_n^\perp}(W)$ we first consider the functions
\begin{eqnarray}
\nonumber
P_{V_n^\perp} \omega_i^*, \quad i=1,\dots, m
\end{eqnarray}
These functions are $0$ for $i=1,\dots,p$ since $\omega_i^*\in V_n$ for these values of $i$. 
They are equal to $\omega_i^*$ for $i=n+1,\dots,m$ and to $\omega_i^*-s_i\phi_i^*$ 
for $i=p+1,\dots,n$, and these $m-p$ functions are non-zero pairwise orthogonal. Therefore an orthonormal basis of $P_{V_n^\perp}(W)$
is given by the normalized functions
\begin{eqnarray}
\nonumber
(1-s_i^2)^{-1/2}(\omega_i^*-s_i\phi_i^*) , \quad i=p+1,\dots,n,\quad {\rm and} \quad \omega_i^*, \quad i=n+1,\dots,m.
\end{eqnarray}
By a similar construction, we find that an orthonormal basis of $P_{W^\perp}(V_n)$ is given by the normalized functions
\begin{eqnarray}
\nonumber
(1-s_i^2)^{-1/2}(\phi_i^*-s_i\omega_i^*),  \quad i=p+1,\dots,n.
\end{eqnarray}
Therefore bases for $V_n^\perp$ and $W^\perp$ are defined as union of these bases with the basis $(\psi_i^*)_{i\geq 1}$
for $V_n^\perp \cap W^\perp$.

Finally, we close out this section, by  giving a parametric description of the set  $\cK_w=\cK_w(V_n)$ for the single space problem  which shows in particular that this set is an ellipsoid.
   \begin{lemma}
  \label{lem:ellipsoids}
  Given a single space $V_n\subset \cH$, the body 
  $$\cK_w:=\cK_w(V_n):=\cK_w^{\one}(V_n):=\{u\in\cK^{\one}(V_n):\ P_Wu=w\}$$
   is a non-degenerate ellipsoid contained in the affine space $\cH_w$.
    \end{lemma}  
  
  \noindent
 {\bf Proof:} Using the favorable bases for $W$ and $W^\perp$, we can write any $u\in \cH_w$ as
 \begin{eqnarray}
\nonumber
 u=\sum_{j=1}^m w_j\omega_j^* +\sum_{j=p+1}^n x_j (1-s_j^2)^{-1/2}(\phi_j^*-s_j\omega_j^*)+\sum_{i\geq 1} y_j \psi_j^*,
 \end{eqnarray}
 where the $w_j=\<w,\omega_j^*\>$ for $j=1,\dots,m$, are given, and the $x_j$ and $y_j$ are the coordinates of $u-w$ in the
 favorable basis of $W^\perp$. We may now write
 $$
 \begin{disarray}{ll}
 P_{V_n^\perp} u &=\sum_{j=1}^m w_jP_{V_n^\perp} \omega_j^*+\sum_{j=p+1}^n x_j (1-s_j^2)^{-1/2}P_{V_n^\perp} (\phi_j^*-s_j\omega_j^*)
 +\sum_{i\geq 1} y_j \psi_j^*\\
  & = \sum_{j=p+1}^m w_j(\omega_j^*-s_j \phi_j^*)-\sum_{j=p+1}^n x_j (1-s_j^2)^{-1/2}s_j(\omega_j^*-s_j \phi_j^*)
 +\sum_{i\geq 1} y_j \psi_j^*\\
 & = \sum_{j=n+1}^m w_j(\omega_j^*-s_j \phi_j^*)+\sum_{j=p+1}^n(w_j- x_j s_j(1-s_j^2)^{-1/2})(\omega_j^*-s_j \phi_j^*)
 +\sum_{i\geq 1} y_j \psi_j^*.
 \end{disarray}
 $$
 All terms in the last sum are pairwise orthogonal and therefore
 \begin{eqnarray}
\nonumber
 \|P_{V_n^\perp} u\|^2 =\sum_{j=n+1}^m (1-s_j^2)w_j^2+\sum_{j=p+1}^n (1-s_j^2)(w_j- x_j s_j(1-s_j^2)^{-1/2})^2+\sum_{j\geq 1} y_j^2.
 \end{eqnarray}
 Now $u\in \cK_w$ if and only if  $\|P_{V_n^\perp} u\|^2 \leq \e_n^2$,  or equivalently
 \be
 \sum_{j=p+1}^n s_j^2(x_j-a_j)^2+\sum_{j\geq 1} y_j^2 \leq  C,
 \label{ellips}
 \ee
 with $C:=\e_n^2-\sum_{j=n+1}^m (1-s_j^2)w_j^2$ and $a_j:=(1-s_j^2)^{1/2}s_j^{-1} w_j$
 which is the equation of a non-degenerate ellipsoid in $\cH_w$.
 \hfill $\Box$
 
 \begin{remark}
 The above equation \iref{ellips} directly shows that the radius of $\cK_w$ is equal to $s_n^{-1}C^{1/2}$
 which is an equivalent expression of \iref{Rstar}.
 \end{remark}

\section{The multi-space problem}
\label{sec:multiple}

  In this section, we consider the multi-space problem  as described in the introduction.   We are interested in the optimal recovery of the elements in the set $\cK:=\cK^{\mult}$ as described
by \eref{defKm}.   For any given $w\in W$, we consider the set %
 \begin{eqnarray}
\nonumber
\cK_w:=\cK_w^{\mult}:=\cK^{\mult}\cap 
{\cH_w}
 =\bigcap_{j=0}^n \cK^j_w,
\end{eqnarray}
where 
\begin{eqnarray}
\nonumber
\cK^j_w:=\cK^j\cap \cH_w:=\{ u\in \cH_w\; : \; \dist(u,V_j)\le \e_j\}.
\end{eqnarray}
 In other words, $\cK^j_w$ is the set in the 
 {one-space}
  problem
considered in the previous section.
We have seen that $\cK^j_w$ is an ellipsoid with known center $u^*_j=u_j^*(w)$
and known Chebyshev radius given by \iref{Rstar} with $n$ replaced by $j$,
{and $u^*$ and $v^*$ replaced by $u_j^*$ and $v_j^*$}
 in that formula.

Thus, $\cK_w $ is now the intersection of $n+1$ ellipsoids.
The optimal algorithm $A^*$, for the recovery of $\cK_w$,  is the one that would find the center of the Chebyshev ball of 
this set and its performance would then be given by its Chebyshev radius.
In contrast to the one-space problem, this center and radius do not have simple
computable expressions. The first results of this section provide an a priori estimate of the Chebyshev radius in the multi-space setting by exploiting favorable bases.   This a priori analysis illustrates when a gain in performance is guaranteed to occur, although the a priori estimates we provide may be pessimistic.

We then give examples which show that the Chebyshev radius in the multi-space case can be   far smaller than the
minimum of the Chebyshev radii of the $\cK^j_w$ for $j=0,\dots,n$. These examples are intended to illustrate that exploiting the multi-space case can be much more advantageous than simply executing the one-space algorithms and taking the one with best performance, see \eref{perfba}.

The latter part of this section proposes two simple algorithmic strategies, each of them converging to a point in $\cK_w$.
These algorithms thus produce a near optimal solution, in the sense that if $A$ is the map
corresponding to either one of them, we have
\be
E_A(\cK_w)\leq 2E_{A^*}(\cK_w)=2E(\cK_w), \quad w\in W,
\label{nearopt1}
\ee
and in particular
\be
E_A(\cK)\leq 2E(\cK).
\label{nearopt2}
\ee
Both of these algorithms are iterative and based on alternating projections.  An a posteriori estimate for the distance between a given
iterate and the intersection of the ellipsoids is given and used both, as a stopping criteria and to analyze the convergence rates of the
algorithms.

\subsection{A priori bounds for  the radius of $\cK_w$}
\label{ss;boundrad}

In this section, we derive a priori bounds for $\rad(\cK_w^{\mult})$. Although these bounds may   overestimate  $\rad(\cK_w^{\mult})$, they allow us to show examples where the multi-space algorithm
is significantly better than simply chosing one space and using the one-space algorithm.  Recall that for the one-space problem,
we observed that $\rad(\cK_w^{\one})$ is largest when $w=0$. The following results 
{show}  that for the multi-space problem
$\rad(\cK_w^{\mult})$ is also controlled by $\rad(\cK_0^{\mult})$, up to a multiplicative constant. Note that $\cK_w^{\mult}$ is generally not a symmetric set,
except for $w=0$.  In going further in this section $\cK$ and $\cK_w$ will refer to the multi-space sets.
\begin{lemma}
\label{lem:rad}
For the multi-space problem, one has
\be
\label{zerobiggest}
\rad(\cK_w)\le 2\rad(\cK_0), \quad w\in W.
\ee
Therefore,
\be
\label{zerobiggest1}
E(\cK)\le 2\rad(\cK_0).
\ee
\end{lemma}

\noindent
{\bf Proof:} 
Fix $w\in W$ and let  $\t u:=\t u(w)$ be the center of the Chebyshev ball for $\cK_w$ which by Remark \ref{rem:center}, belongs to $\cK_w$. For 
any  $u\in \cK_w$ we have $\eta:=\frac{1}{2}(u-\t u)$ is in $W^\perp$ and 
also
$$
\dist (\eta,V_k)\le \frac{1}{2}(\dist(u,V_k)+\dist(\t u,V_k))\le \e_k,\quad k=0,1,\dots,n.
$$
Hence, $\eta\in \cK_0$ which gives
\begin{eqnarray}
\nonumber
\|u-\t u\|=2\|\eta\|\le 2\rad(\cK_0),
\end{eqnarray}
where we have used the fact that, by Remark \ref{rem:diam2r}, the best Chebyshev ball for $\cK_0$ is centered at $0$. This proves
\eref{zerobiggest}.  The estimate \eref{zerobiggest1} follows from the definition of $E(\cK)$.
\hfill $\Box$
\nl

In view of the above Lemma \ref{lem:rad}, we concentrate  on deriving a priori bounds for the radius of the set $\cK_0$.   
We know that $\cK_0$ is the intersection of the ellipsoids $\cK_0^{j}$ for ${j}=0,1,\dots,n$, each of which is centered at zero.   We also know that the Chebyshev ball for $\cK_0^{j}$ is  $B(0,\rad(\cK_0^{j})$ and we know  from \eref{perfba} that 
\begin{eqnarray}
\nonumber
\rad(\cK_0^{j})=  \mu(V_{j},W) \e _{j}, \quad { j}=0,1,\dots,n,
\end{eqnarray}
which is a computable quantity.
This gives the obvious bound
\be
\label{oldestimate}
\rad(\cK_0)\le \min_{0\le  k\le n}  \mu(V_k,W)\e_k.
\ee
In the following, we show  that we can improve on this bound considerably. Since $\cK_0$ is symmetric around the origin, we have 
\begin{eqnarray}
\nonumber
\rad(\cK_0)=\argmax_{\eta\in \cK_0}\|\eta\|.
\end{eqnarray}
So we are interested in bounding $\|\eta\|$ for each $\eta\in \cK_0$.

Since the spaces $V_j$ are nested,
we can consider an orthonormal basis $\{\phi_1,\dots,\phi_n\}$ for $V_n$, for which, 
$\{\phi_1,\dots, \phi_j\}$ is an orthonormal basis for each of the $V_j$ for $j=1,\dots,n$. 
We will use the favorable bases constructed in the previous section 
in the case of the particular space $V_n$. Note that if $\{\phi_1^*,\dots,\phi_n^*\}$ is the favorable basis
for $V_n$, we do not generally have that $\{\phi_1^*,\dots,\phi_j^*\}$ is a basis of $V_j$.

Let  $\eta$ be any element from $\cK_0$.
Since ${\rm dist}(\eta,V_n)\leq \e_{n}$, we may express
$\eta$  as
\begin{eqnarray}
\nonumber
\eta=\sum_{j=1}^n \eta_j \phi_j^* +e=\sum_{j=1}^n \alpha_j \phi_j +e ,\quad e\in V_n^\perp \ {\rm and} \ \|e\|\le \e _n.
\end{eqnarray}
So, 
\begin{eqnarray}
\nonumber
\|\eta\|^2=\sum_{j=1}^n \eta_j^2  +\|e\|^2=\sum_{j=1}^n \alpha_j^2  +\|e\|^2.
\end{eqnarray}
The $\alpha_j$ and $\eta_j$ are related by the equations
\begin{eqnarray}
\nonumber
\sum_{j=1}^n  \lambda_{i,j} \alpha_j=\eta_i, \quad i=1,\dots,n,
\end{eqnarray}
where
\begin{eqnarray}
\nonumber
\lambda_{i,j}:= \langle \phi_j,\phi_i^*\rangle ,\quad 1\le i,j\le n.
\end{eqnarray}
The fact that ${\rm dist}(\eta,V_k)\leq \e_k$ for $k=0,\dots,n$ is expressed by the inequalities
\begin{eqnarray}
\nonumber
\sum_{j=k+1}^n \alpha_j^2+ \|e\|^2\leq \e_k^2, \quad k=0,\dots,n.
\end{eqnarray}
 Since $\eta\in W^\perp$,  we have  that
\begin{eqnarray}
\nonumber
 0=P_W\eta=\sum_{j=1}^n s_j\eta_j  \omega_j^*+P_W e   \,.
 \end{eqnarray}
 It follows that 
\begin{eqnarray}
\nonumber
 \sum_{j=1}^n s_j^2\eta_j^2  = \|P_We\|^2 \leq \|e\|^2 \leq \e_n^2.
 \end{eqnarray}
We now return to the representation of 
 {$\cK_0$}
  in the $\phi_j$ coordinate system.  We know that all $\alpha_j$ satisfy $|\alpha_j|\le \e _{j-1} $.  This means that the coordinates $\{\alpha_1,\dots,\alpha_n\}$ of any point
 in 
 {$\cK_0$}
  are in the $n$-dimensional rectangle 
\begin{eqnarray}
\nonumber
 R=[-\e _0,\e _0]\times \cdots \times [-\e _{n-1},\e _{n-1}].
 \end{eqnarray} 
 It follows that each 
 $\eta_i$
  satisfies the crude estimate
  \be
 \label{betaestimate}
 |\eta_i|\le \sum_{j=1}^n |\lambda_{i,j}||\alpha_j|\le  \sum_{j=1}^n |\lambda_{i,j}|\e _{j-1}=:\theta_i \quad i=1,\dots,n .
 \ee 
The numbers 
$\theta_i$
are computable.  The bound \eref{betaestimate} allows us to estimate  
\begin{eqnarray}
\nonumber
\rad (\cK_0)^2= \sup_{\eta\in\cK_0}\|\eta\| ^2\le \e _n^2+\sup\Big \{ \sum_{j=1}^n\eta_j^2: \ |\eta_j|\le \theta_j\quad {\rm and} \quad \sum_{j=1}^n s_j^2\eta_j^2 \le \e _n^2\Big\}
 \end{eqnarray}
 Since the  $s_{j}$ are non-increasing,  the supremum on the right side takes the form
\begin{eqnarray}
\nonumber
 \delta\theta_{k}^2+ \sum_{j=k+1}^n\theta_{j}^2,\quad 0<\delta\le 1,
 \end{eqnarray}
 where $k$ is   the largest integer such  that
  \be
 \label{right1}
 \sum_{j=k}^ns_j^2\theta_{j}^2\ge \e _n^2,
 \ee
 and $\delta$ is chosen so that  
  \be
 \label{right2}
\delta s_k^2\theta_k^2  + \sum_{j=k+1}^n s_j^2\theta_{j}^2= \e _n^2.
 \ee
This gives us the following  bound on the Chebyshev radius of $\cK_0$.
\be
\label{bounddiameter}
\rad(\cK_0)^2\le \e _n^2+ \delta \theta^2_{k}+  \sum_{j=k+1}^n  \theta_{j}^2 :=E_n^2.
 \ee
 Using this estimate together with Lemma \ref{lem:rad}, we have proven the following theorem.

\begin{theorem}
\label{thm:ms}
For the multi-space problem, we have the following estimates for Chebyshev radii.  For $\cK_0$, we have
\begin{eqnarray}
\nonumber
\rad(\cK_0)\le E_n,
\end{eqnarray}
where $E_n:=\Big(\e _n^2+ \delta \theta^2_{k}+  \sum_{j=k+1}^n  \theta_{j}^2\Big)^{1/2}$.   For any $w\in W$, we have
\begin{eqnarray}
\nonumber
\rad(\cK_w)\le 2 E_n.
\end{eqnarray}
For $\cK$, we have the bound
\begin{eqnarray}
\nonumber
\rad(\cK)\le 2 E_n.
\end{eqnarray}
\end{theorem}

 We next compare the bound in \eref{bounddiameter} with the one space bound %
\begin{eqnarray}
\nonumber
\rad(\cK_{0})\le\mu(V_n,W)\e _n= s_n^{-1}\e _n ,
\end{eqnarray}
 which is obtained by considering only the approximation property of $V_n$ and not exploiting the other spaces $V_j$, $j<n$, see \eref{oldestimate}.  For this, we return to 
 the definition of $k$ from \eref{right1}.    We can write each term that appears in \eref{right2} as $\gamma_j\e _n^2$ where 
 {$\sum_{j=k}^n\gamma_j=1$}.  In other words,
 \begin{eqnarray}
\nonumber
 \theta_{j}^2=\gamma_js_j^{-2}\e _n^2, \quad k<j\le n,\quad \theta_{k}^2=\delta^{-1}\gamma_{k}s_k^{-2}\e _n^2.
 \end{eqnarray}
 Hence,
\begin{eqnarray}
\nonumber
 E_n^2\le \e _n^2+s_n^{-2}\e _n^2 \leq 2s_n^{-2} \e_n^2,
 \end{eqnarray}
which is at least as good as the old bound up to a multiplicative constant $\sqrt 2$. 

We finally observe that the bound $E_n$ is obtained by using 
the entire sequence $\{V_0,\dots,V_n\}$. Similar bounds $E_\Gamma$ are obtained
when using a subsequence $\{V_j \; : \; j\in \Gamma\}$ for any $\Gamma \subset \{0,\dots,n\}$.
This leads to the improved bound
\begin{eqnarray}
\nonumber
\rad(\cK_0)\le 
{\min} \{E_\Gamma\; : \; \Gamma\subset \{0,\dots,n\}\}.
\end{eqnarray}
In particular defining $E_j=E_\Gamma$ for $\Gamma=\{0,\dots,j\}$ we find that 
\begin{eqnarray}
\nonumber
 E_j^2\le 2\mu(V_j,W)^2\e _j^2.
 \end{eqnarray}
 Therefore
\begin{eqnarray}
\nonumber
 E_n^* \leq \sqrt 2 \min_{j=0,\dots,n} \mu(V_j,W)\e _j,
 \end{eqnarray}
 which shows that the new estimate is as good as \iref{oldestimate} up to the multiplicative
 constant $\sqrt 2$.
  
\subsection{Examples}
\label{ss:examples}

One can easily find examples for which   the Chebyshev radius of $\cK_w$ is substantially smaller than
the minimum of the Chebyshev radii of the $\cK^j_w$, therefore giving higher potential accuracy 
in the multi-space approach. As a simple example to begin this discussion,  consider the case where 
\begin{eqnarray}
\nonumber
\cH=\R^2 , \quad V_0=\{0\}, \quad V_1=\R e_1, \quad W=\R (e_1+e_2)
\end{eqnarray}
where $e_1=(1,0)$ and $e_2=(0,1)$. So, $V_1$ and $W$ are one dimensional spaces. Then, with the choices
\begin{eqnarray}
\nonumber
\e_0=1,\quad \e_1=\frac 1 2, \quad w=\(\frac {\sqrt 3 +1} 4,\frac {\sqrt 3 +1} 4\) ,
\end{eqnarray}
it is easily seen that $\cK_w$ is the single point $\(\frac {\sqrt 3} 2, \frac 1 2\)$ and has therefore
null Chebyshev radius while $\cK^0_w$ and $\cK^1_w$ have positive Chebyshev radii.

In more general settings we do not have such a simple description of $\cK_w$, however we
now give some additional examples that show that even the a priori estimates of the previous 
section can be significantly better than the one space estimate as well as the estimate \eref{oldestimate}.  We  consider the two
 extremes in the compatibility between the favorable basis $\{\phi_1^*,\dots,\phi_n^*\}$ and the basis 
 $\{\phi_1,\dots,\phi_n\}$  which describes the approximation properties of the sequence $\{V_0,\dots,V_n\}$.    
 \nl
 \nl
  {\bf Example 1:}  In this example we consider the case where the two bases coincide,
\begin{eqnarray}
\nonumber
 \phi_i^*=\phi_i,  \quad i=1,\dots,n.  
 \end{eqnarray}
 Note that in this case the singular values $\{s_1,\dots,s_k\}$ for the pair $\{V_k,W\}$ coincide with the 
 first $k$ singular values for the pair $\{V_n,W\}$. Therefore
\begin{eqnarray}
\nonumber
\mu(V_k,W)=s_k^{-1}, \quad k=0,\dots,n,
 \end{eqnarray}
 where we have set $s_0:=1$. We also have
\begin{eqnarray}
\nonumber
 \theta_k=\e _{k-1},\quad k=1,\dots,n.
 \end{eqnarray}
  We fix $\e_n:=\e$ and $\e_{n-1}:=\e_{n-2}:=\e^{1/2}$  and the values $s_n:= \e$ and $s_{n-1}:=s_{n-2}:=\e^{1/2}$ and 
 all other $\e_k:=1$ and all other $s_k:=1$.  We examine what happens when $\e$ is very small.   The  {estimate}  \eref{olderror}   would give the bound 
 \begin{eqnarray}
\nonumber
\min_{0\le k\le n} \mu(V_k,W) \e_k=\min_{0\le k\le n} s_k^{-1}\e_k =1,
 \end{eqnarray}
 as the bound for $\rad(\cK_0)$ and $E(\cK)$.
On the other hand, since, 
\begin{eqnarray}
\nonumber
s_{n}^{2}\e_{n-1}^2=\e^3\ll \e^2 \quad {\rm and}\quad  s_{n-1}^{2}\e_{n-2}^2= \e^2,
\end{eqnarray}
the value of $k$ in \eref{right1} is $n-1$.  It follows that the error $E_n$ in the multi-space method   \eref{bounddiameter} satisfies
\begin{eqnarray}
\nonumber
E_n^2\le  \e_{n-2}^2+ \e_{n-1}^2+ \e_n^2\le 3\e.
\end{eqnarray}
Hence, the error for the multi-space method can be arbitrarily small as compared to the error of the one-space method.
 \nl
 \nl
 {\bf Example 2: } 
  We next consider the other extreme where the two bases are incoherent in the sense that each entry in the change of basis matrix
  satisfies
 \begin{eqnarray}
\nonumber
  |\lambda_{i,j}|\le C_0n^{-1/2},\quad 1\le i,j\le n.
  \end{eqnarray}
  We want to show that $E_n$ can be smaller 
  than
  the estimate in \iref{oldestimate} in this case as well.
  To illustrate how the estimates go, we assume that $n\ge 2$ and  $|\lambda_{i,j}|=1/\sqrt{n}$,  for all $1\le i,j\le n$.    We  will
 take
 \begin{eqnarray}
\nonumber
s_n\ll  s=s_1=s_2=\dots=s_{n-1},
 \end{eqnarray}
with the values of $s$ and $s_n$ specified below.  We define
\begin{eqnarray}
\nonumber
\e_0:=1/2\quad {\rm and}\quad  \e_j=\frac {1}{2(n-1)},\quad j=1,\dots,n-1,
 \end{eqnarray}
 so that $\sum_{j=0}^{n-1}\e_j=1$.   
 It follows from the definiton of $\theta_k$ given in \eref{betaestimate} that
    
\begin{eqnarray}
\nonumber
    \theta_k=1/\sqrt{n}:=\theta, \quad  k=1,\dots,n.
\end{eqnarray}

With these choices, the best one space estimate   \eref{olderror} 
{is}
   \be
   \label{bestone}
   \min\{\e_0,s^{-1}\e _{n-1},s_n^{-1}\e_n\}. 
   \ee
 Now, we take $\e_n$ very 
 {small}
  and  $s_n=\e_n^2$.  We then choose $s$ so that 
  \be
  \label{case2}
(s^2+s_n^2) \theta^2 =\e_n^2.
  \ee
  This gives  $k=n-1$ in \eref{right1}    and so 
 \begin{eqnarray}
\nonumber
  E_n^2=\e_n^2+\theta_{n-1}^2+\theta_n^2 \leq 3n^{-1}.
  \end{eqnarray}
 On the other hand, \eref{case2} says that 
 {$s^{-1}=\e_n^{-1}(n-\e_n^2)^{-1/2}$}.
 Thus,  from \eref{bestone}, the best one space estimate 
 {is} 
  \begin{eqnarray}
\nonumber
   \min\{\e_0,s^{-1}\e _{n-1},s_n^{-1}\e_n\} 
   = 
   \min \Big\{\frac{1}{2}, 
   {\frac{1}{2(n-1)\sqrt{n-\e_n^2}}}
   \e_n^{-1},  \e_n^{-1} \Big \} =1/2, 
   \end{eqnarray}
provided {$\e_n\leq  n^{-3/2}$.}
 Hence, the multi-space estimate   \eref{bounddiameter} is better than the one space estimate by at least the
factor $n^{-1/2}$ in this case.

 \subsection{Numerical algorithms} 
 
 \label{ss:numerical}   
 
 In this section, we discuss  some possible numerical algorithms, based on convex optimization, 
 for the multi-space case. For any given 
 data $w\in W$,
 such that $\cK_w$ is not empty, these algorithms produce, in the limit,
 an element $A(w)$ which belongs to $\cK_w$, so that they are
 near optimal in the sense of \iref{nearopt1} and \iref{nearopt2}.
 
 We recall that $\cK_w$ is given by
\begin{eqnarray}
\nonumber
 \cK_w=\cH_w\cap \cK^0 \cap \cK^1\cap \cdots \cap \cK^n.
 \end{eqnarray}
 One first observation is that although the set $\cK_w$ 
 {may be}
 infinite dimensional, we may reduce the
 search for an element in $\cK_w$ to the finite dimensional  space
\begin{eqnarray}
\nonumber
 \cF:=V_n+ W,
 \end{eqnarray}
 which has dimension $d=m+n-p$, where $p=\dim(V_n\cap W)$. Indeed, if $u\in \cK_w$, then its projection
 $P_\cF u$ onto $\cF$ remains in $\cK_w$, since $u-P_\cF u\in W^\perp\cap V_n^\perp$ implies
\begin{eqnarray}
\nonumber
 P_W P_\cF u=P_Wu=w,
 \end{eqnarray}
 and
\begin{eqnarray}
\nonumber
 \dist(P_\cF u, V_j)\leq  \dist( u, V_j)\leq \e_j, \quad j=0,\dots,n.
 \end{eqnarray}
 Therefore, without loss of generality, we may assume that 
\begin{eqnarray}
\nonumber
 \cH=\cF,
 \end{eqnarray}
and that  the sets $\cH_w$ and {$\cK^j$} that define $\cK_w$ are contained in this finite dimensional space.
 
The problem of finding a point in the intersection of convex sets is sometimes referred to as {\it convex feasibility}
and has been widely studied in various contexts. We refer to \cite{Comb1,Comb2} for surveys on
various possible algorithmic methods. We restrict our
discussion to two of them which have very simple expressions in our particular case. Both are based 
on the orthogonal projection operators onto the spaces $\cH_w$ and $\cK^j$.  Let us first observe that these projections
are very simple to compute. For the projection onto $\cH_w$, we use the orthonormal basis $\{\omega_1,\dots,\omega_m\}$ of $W$.
For any $u\in \cF$ we have
\be
P_{\cH_w} u= P_{W^\perp}u+w=u-\sum_{i=1}^{m} \<u,\omega_i\>\omega_i+w.
\label{PHw}
\ee
For the projection onto $\cK^j$, we extend the basis $\{\phi_1,\dots,\phi_n\}$ into  
an orthonormal basis $\{\phi_1,\dots,\phi_d\}$ of $\cF$.  We then have
\begin{eqnarray}
\nonumber
P_{\cK^j} u= \sum_{i=1}^j \<u,\phi_i\> \phi_i+\alpha\(\sum_{i=j+1}^d  \<u,\phi_i\> \phi_i\), \quad \alpha:=\min\Big\{1,\e_j\(\sum_{i=j+1}^{d}Ê|\<u,\phi_i\>|^2\)^{-1/2}\Big\}.
\end{eqnarray}
We 
{now}
 describe two elementary and well-known algorithms.
\nl
\nl
{\bf Algorithm 1: sequential projections.}  This algorithm  is  a cyclical application of the above operators.   Namely, starting
say from $u^0=w$, we define for $k\geq 0$ the iterates
\begin{eqnarray}
\nonumber
u^{k+1}:=P_{\cK^n}P_{\cK^{n-1}} \cdots P_{\cK^1}P_{\cK^0} P_{\cH_w} u^{k}.
\end{eqnarray}
We know from general results on alternate projections onto convex sets \cite{Br} that this sequence converges towards
a point $u^*\in \cK_w$ when $\cK_w$ is not empty. We make further use of the
following observation: the nestedness 
property $V_0\subset V_1 \subset \dots \subset V_n$ implies 
that $u^k$ belongs to $\cK=\cK^0\cap \dots \cap \cK^n$.
\nl
\nl
{\bf Algorithm 2: parallel projections.} This algorithm combines the projections onto the sets $\cK$ according to 
\begin{eqnarray}
\nonumber
u^{k+1}:=P_{\cH_w}\(\sum_{j=0}^n \gamma_jP_{\cK^j}\) u^{k},
\end{eqnarray}
where the weights $0<\gamma_j<1$ are such that $\gamma_0+\cdots +\gamma_n=1$, for example 
$\gamma_j:=\frac 1 {n+1}$. It may be
viewed as a projected gradient iteration for the minimization over $\cH_w$ of the differentiable function
\begin{eqnarray}
\nonumber
F(u):=\sum_{j=0}^n \gamma_jF_j(u),\quad F_j{(u)}:=\frac 1 2\dist(u,\cK^j)^2.
\end{eqnarray}
Notice that the minimum of $F$  is attained exactly at each point of $\cK$. Since $\nabla F_j(u)=u-P_{\cK^j}u$, we find that
\begin{eqnarray}
\nonumber
u^{k+1}=P_{\cH_w}(u^k-\nabla F(u^k)).
\end{eqnarray}
Classical results on constrained minimization methods \cite{LP} show that this algorithm converges toward 
a minimizer $u^*$ of $F(u)$ over $\cH_w$ which
clearly belongs to $\cK_w$ when $\cK_w$ is not empty.

\subsection{A posteriori estimate and convergence rates}

Each of the above algorithms generates a sequence $(u^k)_{k\ge 1}$ of elements from $ \cF$ which are guaranteed to converge to a point
in $\cK_w$ provided that  this set is nonempty.   We would like
to
 have a bound for $\dist(u^k,\cK_w)$, since this would allow us to check
the progress of the algorithm and also could be  utilized 
as 
a stopping criterion when we have gained sufficient accuracy.
Here we restrict our analysis to Algorithm 1.

We will use certain geometric properties of the set $\cK$, expressed by the following
lemma.
\begin{lemma}
\label{lem:ecc}
If $u_1,u_2\in\cK$  then the ball $B:=B(u_0,r)$ centered at $u_0:=\frac{1}{2}(u_1+u_2)$ of radius   
\be
\label{radius}
r:=\frac{1}{8}\min_{j=0,\dots,n} \e_j^{-1}\|P_{V_j^\perp} (u_1)- P_{V_j^\perp} (u_2)\|^2     
\ee
is completely contained in $\cK$.
\end{lemma}
\noindent
{\bf Proof:}  For  $u_1,u_2\in \cK^j$ the ball $B(u_0,r)$ is contained in $\cK^j$ if and only if the ball in $V_j^\perp$  centered at $P_{V_j^\perp}u_0$ with the radius $r$ is contained in 
$P_{V_j^\perp}(\cK^j)=\{x\in V_j^\perp \ :\ \|x\|\leq \epsilon_j\}:=\cB_j$. Let $v^j_s:=P_{V_j^\perp}(u_s)$ for $s=0,1,2$ and let  $\delta_j:=\|v_1^j-v_2^j\|$. The parallelogram identity gives
$$
\|v_0^j\|^2=\frac12\| v_1^j\|^2+\frac12\|v_2^j\|^2 -\frac14\|v_1^j-v_2^j\|^2,
$$
so that $\|v_0^j\|^2\leq \e_j^2-\frac14 \delta_j^2$.  Thus for 
\begin{eqnarray}
\nonumber
r_j:=\e_j -\sqrt{\e_j^2-\frac14 \delta_j^2}=\e_j\left(1-\sqrt{1-\frac {\delta_j^2}{4\e_j^2}}\right),
\end{eqnarray}
the ball in $V_j^\perp$ centered at $v_0^j$ with radius $r_j$ is contained in $\cB_j$.  Thus, with 
\begin{eqnarray}
\nonumber
\rho:=\min_{j=0,1,\dots,n} r_j,
\end{eqnarray}
we have $B(u_0,\rho)\subset \cK$. Since $\delta_j\leq 2 \e_j$ and $(1-\sqrt{1-x})\geq x/2$ for $0\leq x\leq 1$ we get 
$r_j\geq \delta_j^2/(8\e_j)$ and therefore $\rho\geq r$ from which \iref{radius} follows.  \hfill $\Box$
\nl

We have noticed that the iterates $u^k$ of Algorithm 1
all belong to $\cK$ and we would like to estimate their distance 
from the convex set $\cK_w$.    Let $P_{\cK_w}(x)$ denote the point from $\cK_w$ closest to $x$. This is a well defined map. The following result gives
an estimate for the distance of any $u\in \cK$ from $\cK_w$, in terms of its
distance from the affine space $\cH_w$. This latter quantity is
easily computed using \iref{PHw} which shows that
\begin{eqnarray}
\nonumber
u-P_{\cH_w} u=P_Wu-w=\sum_{i=1}^{m} \<u,\omega_i\>\omega_i-w.
\end{eqnarray}

\noindent

\begin{lemma} 
\label{lem:apost}
Let $u\in \cK$ be such that
\begin{eqnarray}
\nonumber
\alpha:=\dist(u,\cH_w)>0.
\end{eqnarray}
Then
\be
\label{toprove}
\| P_{\cH_w} u-P_{\cK_w} u\| \leq \rho=\rho(\alpha):=  \max_j   \mu_j(\alpha +4\sqrt{\alpha \e_j}),
\ee
where $\mu_j=\mu(V_j, W)$.
Since $u-P_{\cH_w}u$ is orthogonal to $P_{\cH_w}u-P_{\cK_w} u$, we have
\begin{eqnarray}
\nonumber
\dist(u,\cK_w)^2 \leq \alpha^2+\rho(\alpha){^2}.
\end{eqnarray}
\end{lemma}

\noindent
{\bf Proof:} We set  $u_2=P_{\cK_w} u$ and $\eta=u-u_2$ which we decompose as
\begin{eqnarray}
\nonumber
\eta=(u-P_{\cH_w} u)+(P_{\cH_w} u-u_2)=:\eta_1+\eta_2.
\end{eqnarray}
We wish to show that $\|\eta_2\|\le \rho$, where $\rho$ is defined in \iref{toprove}. To this end, 
observe that  $\eta_1\in W$ and $\eta_2\in W^\perp$ so  that this is an orthogonal decomposition. 
Moreover, using \iref{mu} and noting that  $\|\eta_1\|=\alpha$, we have   
\be
 \label{PW5}
\|P_{V_j^\perp} \eta\|\geq \|P_{V_j^\perp} \eta_2\|-\|P_{V_j^\perp}\eta_1\|\geq \beta(V_j,W)\|\eta_2\|-\alpha.
\ee
We infer 
from   Lemma \ref{lem:ecc}   that the ball $B$ with center at $u_0=\frac12(u+u_2)$ and radius 
$$
r=\frac18\min_{j=0,1,\dots,n} \e_j^{-1} \|P_{V_j^\perp} \eta\|^2
$$
 is contained in  $\cK$. Let us suppose now that  $\|\eta_2\|>\rho$ and derive a contradiction. 
 Then, we obtain from \eref{PW5}
 \begin{eqnarray}
\nonumber
\|P_{V_j^\perp} \eta\| >  \mu_j^{-1}\rho-\alpha
\geq \mu_j ^{-1}\mu_j(\alpha +4\sqrt{\alpha\e_j})-\alpha=4\sqrt{\alpha \e_j},
\end{eqnarray}
and thus
\begin{eqnarray}
\nonumber
r > \frac18\min_{j=0,1,\dots,n} \e_j^{-1} 16\alpha\e_j =2\alpha.
\end{eqnarray}
On the other hand, 
note that $\|u_0-P_{\cH_w} u_0\|=\frac12 \|u-P_{\cH_w }u\|=\alpha/2$. Therefore, $P_{\cH_w}u_0\in \cK$ and hence in $\cK_w$.
Moreover,
\begin{eqnarray}
\nonumber
\|u- P_{\cH_w}u_0\|^2 = \alpha^2 + \frac 14 \|u_2 - P_{\cH_w}u\|^2,
\end{eqnarray}
and 
\begin{eqnarray}
\nonumber
\|u- u_2\|^2=\alpha^2 + \|u_2 - P_{\cH_w}u\|^2.
\end{eqnarray}
If $u_2\neq P_{\cH_w}u$, we have $\|u- P_{\cH_w}u_0\| < \|u- u_2\|$ which is a contradiction since $u_2$ is the closest point 
to $u$ in $\cK_w$. If $u_2 - P_{\cH_w}u =0$ then $\eta_2=0$ contradicting $\|\eta_2\| >\rho$. This completes the proof.
\hfill $\Box$\\

One 
immediate
consequence of the above lemma is an a posteriori error estimate for 
the squared distance to $\cK_w$ 
\begin{eqnarray}
\nonumber
\delta_k:=\dist(u^k,\cK_w)^2,
\end{eqnarray}
in Algorithm 1. Indeed, we have observed that
$u^k\in \cK$ and therefore
\begin{eqnarray}
\nonumber
\delta_k \leq \alpha^2_k+\rho(\alpha_k)^2, \quad \alpha_k:=\dist(u^k,\cH_w).
\end{eqnarray}
This ensures  the following  accuracy with respect to the unknown $u\in \cK_w$:
\begin{eqnarray}
\nonumber
\|u-u^k\| \leq \sqrt {\alpha^2_k+\rho(\alpha_k)^2}+2 \rad(\cK_w).
\end{eqnarray}
If we have an a priori estimate for the Chebyshev radius of $\cK_w$, such as the bound $E_n$ from Theorem \ref{thm:ms},
one possible stopping criterion is the validity of
\begin{eqnarray}
\nonumber
 \sqrt {\alpha^2_k+\rho(\alpha_k)^2}\leq E_n.
\end{eqnarray}
This ensures that we have achieved accuracy  $\|u-u^k\|\leq 3E_n$, however note that $E_n$ can sometimes be a very pessimistic
bound for $\rad(\cK_w)$ so that significantly higher accuracy is reachable by more iterations.

We can also use Lemma \ref{lem:apost} to establish a convergence estimate for
$\delta_k$ in Algorithm 1. For this purpose, we introduce the
intermediate iterates
\begin{eqnarray}
\nonumber
u^{k+\frac 1 2}:= P_{\cH_w} u^{k},
\end{eqnarray}
and the corresponding squared distance 
\begin{eqnarray}
\nonumber
\delta_{k+\frac 1 2}{:=}\dist(u^{k+\frac 1 2},\cK_w)^2.
\end{eqnarray}
Since the distance to $\cK_w$ is non-increasing 
in each projection steps, it follows that 
\begin{eqnarray}
\nonumber
\delta_{k+1}\leq \delta_{k+\frac 1 2}= \delta_k-\alpha_k^2.
\end{eqnarray}
On the other hand, it easily follows  from Lemma \ref{lem:apost} that 
\begin{eqnarray}
\nonumber
\delta_k-\alpha_k^2 \leq \rho(\alpha_k)^2\leq  A\alpha_k,
\end{eqnarray}
where $A$ is a constant depending on $\epsilon_j$'s, $\mu_j$'s and $\|u\|$. It is easily seen that this implies the validity of the  inequality
\begin{eqnarray}
\nonumber
\alpha_k\geq\sqrt{\delta_k+A^2/4}-A/2\geq \frac{\delta_k}{\sqrt{A^2+4\delta_k}}\geq  \frac{\delta_k}{\sqrt{A^2+4\delta_0}}:=c\delta_k,
\end{eqnarray}
and {therefore}
\begin{eqnarray}
\nonumber
\delta_{k+1} \leq \delta_k-c^2 \delta_k^2.
\end{eqnarray}
From this, one finds by induction that 
\begin{eqnarray}
\nonumber
\delta_k \leq Ck^{-1}, \quad k\geq 1,
\end{eqnarray}
for a suitably chosen constant 
$C:=\max\{c^{-2},\delta_1\}$ taking into account that for any $t\ge1$
\begin{eqnarray}
\nonumber
\frac{C}{t}\(1-\frac{Cc^2}{t}\)\le C\(\frac{t-1}{t^2}\) \le \frac{C}{t+1} \,.
\end{eqnarray}
\hfill $\Box$

\begin{remark}
The above convergence rate $\cO(k^{-1/2})$ for the distance between $u^k$ and $\cK_w$ is quite pessimistic,
however, one can easily exhibit examples in which it 
indeed occurs due to the fact that $\cH_w$ intersects $\cK$ at a single point
of tangency. On the other hand, one can also easily find other examples for which
convergence of Algorithm 1 is exponential. In particular, this occurs 
whenever $\cK_w$ has an element lying in the interior of $\cK$.
\end{remark}

\end{document}